\documentclass[12pt,leqno]{article}
\usepackage{amsmath,amsfonts,amscd,ProtocoloTope,url}
\usepackage[applemac]{inputenc}

\title{\Large\bf  Cohomological tautness  for Riemannian foliations}

\author{
José Ignacio Royo Prieto\thanks{Departamento de Matemática
Aplicada UPV-EHU Alameda de Urquijo s/n
48013 Bilbao Spain.  {\sl joseignacio.royo@ehu.es}.  Partially
supported by the UPV-EHU grant EHU06/05. }
\\
{\small Universidad del País Vasco} \\ [-0,2cm] {\small Euskal
Herriko Unibertsitatea} \and Martintxo
Saralegi-Aranguren\thanks{Fédération CNRS Nord-Pas-de-Calais FR
2956. UPRES-EA 2462 LML. Faculté Jean Perrin. Université d'Artois.
Rue Jean Souvraz SP 18.   62 307 Lens Cedex - France. {\sl
saralegi@euler.univ-artois.fr}.  }
\\ {\small Université d'Artois }
\and
Robert Wolak\thanks{Instytut Matematyki. Uniwersytet Jagiellonski.
Wladyslawa Reymonta 4, 30 059 Krakow - Poland.
{\sl robert.wolak@im.uj.edu.pl} Partially supported by the KBN grant 2 PO3A
021 25.}
\\ {\small Uniwersytet Jagiellonski}
}

\begin{document}  

\maketitle

\noindent To Nicolae Teleman on his 65th birthday.

\begin{abstract}
In this paper we present some new results on the tautness of
Riemannian foliations in their historical context. The first part
of the paper gives a short history of the problem. For a closed manifold, 
the tautness of a Riemannian foliation can be characterized cohomologically. We extend 
this cohomological characterization to a class of foliations which includes the foliated 
strata of any singular Riemannian foliation of a closed manifold. 

    \end{abstract}

Y. Carrière in his paper, cf. \cite{CA}, based on his Ph.D.  thesis  conjectured  that for Riemannian foliations of compact manifolds, the property ``taut'' understood as the existence of a Riemannian bundle-like metric making all leaves minimal is equivalent to the non-triviality of the top dimension basic cohomology group. The conjecture was based on the previous results of A. Haefliger, cf. \cite{HAE_JDG} , demonstrating that ``being taut'' is a transverse property and on his own research into Riemannian flows on 3-manifolds. For over a decade the conjecture was the subject of intensive study 
by a group of ``feuilleteurs'', being finally solved by X. Masa, cf. \cite{MA}, and refined by J.A. Álvarez, \cite{Suso}. 
 The best account of the development of the theory up to 1995 can
be found in Ph. Tondeur's book, cf. \cite{To}.

The case of  non-compact manifolds is much more complicated as the tautness class of a Riemannian foliation cannot be defined in the standard way as in the case of closed manifolds, cf. \cite{CE}. 
However, for some non-compact manifolds it is possible to propose a similar characterization. In a previous paper we proved that if a foliated Riemannian manifold $(M,g,\F)$ can be embedded as a regular stratum of a singular Riemannian foliation (SRF), then  the following conditions are equivalent:

\begin{itemize}

\item[1)] $\mathcal F$ is taut;

\item[2)] $\kappab = 0,$ where $\kappab = [\ib{\kappa}{\mu} ] \in \lau{H}{1}{}{\mf},$ and $\ib{\kappa}{\mu} $ is the mean curvature form of the bundle-like Riemannian metric $\mu$;

\item[3)] $\lau{H}{0}{\ib{\kappa}{\mu} }{\mf} \neq 0$, where $\mu$ is a bundle-like Riemannian metric;

\item[4)]  $\lau{H}{n}{c}{\mf}\ne 0$, where $n =\codim {\mathcal F}$ and the foliation is transversally oriented.

\end{itemize}

In this paper we extend this characterization to a class of non-compact foliated Riemannian manifolds which include
 not only regular strata of SRFs, but other strata as well (cf. Theorem \ref{T1}, \ref{T2} and \ref{T3}).

    \bigskip

    In the sequel $M$ and $N$  are connected, second countable, Haussdorff,
    without boundary and smooth
    (of class $C^\infty$) manifolds  of dimension   $m$.
    All the maps are considered smooth unless something else is indicated.
We consider on $M$ a Riemannian foliation\footnote{For the notions related
    with Riemannian foliations we refer the reader to \cite{Mo,To}.} $\F$ whose codimension is
    $n$. If $V$ is a saturated submanifold of $M$ we shall denote by
    $(V,\F)$ the induced foliated manifold and $\F_{V}$ the induced
    Riemannian foliation.

\smallskip

\section{An historical overview of the problem}

An involutive subbundle $E$ of
dimension $p$ of $TM$ is called a foliation of dimension $p$ and
codimension $n=m-p.$  The foliation $\F$ is said to be
modelled on a $n$-manifold $N_{0}$ if it is defined by a cocycle
${\cal U} = \{ U_{i},f_{i},g_{ij}\}_{I} $ modelled on $N_{0}$,
i.e.

\begin{enumerate}
\item     $\{ U_{i} \}$ is an open covering of $M$,

\item $f_{i} \colon U_{i} \longrightarrow N_{0}$ are submersions
with connected fibres, and

\item  $g_{ij} f_{j} = f_{i}$ on $U_{i} \cap U_{j}$.
\end{enumerate}

The $n$-manifold $T = \coprod T_{i}$, $T_{i} = f_{i}(U_{i})$, is
called the transverse manifold associated to the cocycle $\cal U$
and the pseudogroup $\cal H$ of local diffeomorphisms of $T$
generated by $g_{ij}$ the holonomy pseudogroup representative on
$T$ (associated to the cocycle $\cal U$). $T$ is a complete
transverse manifold. The equivalence class of $\cal H$ we call the
holonomy pseudogroup of $\F$ (or $(M,\F)$). 
It is not difficult to check that  different cocycles defining the
same foliation provide us with equivalent holonomy pseudogroups,
cf. \cite{H1,H2}. In general, the converse is not true.
The notion of a Riemannian foliation was introduced by Bruce
Reinhart in \cite{RE_AM,RE_AJM}.

A foliation $\mathcal F$ on the smooth manifold $M$
is {\em Riemannian} if on $M$ there exists a {\em bundle-like}
metric $\mu$ for the foliation $\mathcal F, $ (i.e., a geodesic
perpendicular to a leaf of $\mathcal F$ at a point remains
perpendicular to every leaf it meets). In a local adapted chart
$(x_1,\dots ,x_{p},y_1,\dots ,y_n)$ the bundle-like metric $\mu$ has
a representation

$$ \suma{ij=1}{p}\mu_{ij}(x,y) v_i\otimes v_j +
\suma{\alpha \beta =1}{n}\mu_{\alpha \beta }(y)dy_{\alpha} \otimes
dy_{\beta} $$

\noindent where $v_i$ is a 1-form annihilating the bundle
$T{\mathcal F}^{\bot}$ and $v_i( \partial/\partial x_j ) =
\delta^i_j.$

%\begin{minipage}{0.49\textwidth}
%\includegraphics{geodesic.EPS}
%\end{minipage}

\medskip

Let $(M,{\mathcal F})$ be a Riemannian foliation with a
bundle-like metric $\mu$. Then it is defined by a cocycle ${\mathcal
U} = \{ U_i,f_i,k_{ij}\}_I$ modelled on a Riemannian manifold
$(N_0,\bar{g})$ such that

\begin{itemize}

\item[(i)] $f_i\colon (U_i,\mu) \rightarrow (T_0,\bar{\mu})$ is a Riemannian
submersion with connected fibres;

\item[(ii)]  $k_{ij}\colon f_j(U_i\cap U_j) \rightarrow f_i(U_i\cap U_j)$
are local isometries of $(T_0,\bar{\mu})$;

\item[(iii)] $f_j=k_{ji}f_i$ on $f_i(U_i\cap U_j)$.
\end{itemize}

\medskip

A foliation $\mathcal F$ on a Riemannian manifold $(M,\mu)$ is said
to be {\it minimal\/} if all its leaves are minimal submanifolds
of $(M,\mu).$ A foliation $\mathcal F$ on  a manifold $M$ is said to
be {\it taut\/} if there exists a Riemannian metric $\mu$ on the
manifold $M$ for which all  leaves are minimal submanifolds of
$(M,\mu).$

Among other things B. Reinhart introduced and studied the basic
cohomology of these foliations.

\medskip

In the presence of the Riemannian metric $\mu$, the tangent bundle
$TM$ admits an orthogonal splitting $TM= T{\mathcal F} \oplus
T{\mathcal F}^{\bot}.$
We say that the  $k$-form $\alpha$ is of pure type
$(r,s),$ r+s=k, if for any point of $M$ there exists an adapted
chart $(x_1,\dots ,x_{p},y_1,\dots ,y_n)$  such that
$$
\alpha = \suma{}{} f_{IJ} v_{i_1}\wedge \cdots \wedge v_{i_r} \wedge
d{y}_{j_1}\wedge \dots \wedge d{y}_{j_s} \ .
$$
where $1 \leq i_1 <\dots <i_r \leq p$, $1\leq j_1
< \dots <j_s \leq n$, $I=(i_1,\dots ,i_r)$, $J=(j_1,\dots ,j_s)$.

Let us denote by $\hiru{\Om}{k}{M}$ the space of
 $k$-forms  $M$, and by $\hiru{\Om}{r,s}{M}$
the space of forms of pure type $(r,s)$. Then
$$
 \hiru{\Om}{k}{M} = \bigoplus \limits_{r+s=k} \hiru{\Om}{r,s}{M} \ ,
$$
 for short $\bi{\Om}{k}= {\displaystyle \bigoplus_{r+s=k}}\bi{\Om}{r,s}$.

The exterior differential $d \colon \hiru{\Om}{k}{M} \longrightarrow
\hiru{\Om}{k+1}{M} $ decomposes itself into three components $d=
d_{\mathcal F} + d_T + \delta$, where $d_{\mathcal F} $ is of
bidegree $(1,0)$ and $d_T$ is of bidegree $(0,1)$, and $\delta$ is
of bidegree $(-1,2),$ i.e.
$$
 d_{\mathcal F} \colon \bi{\Om}{r,s} \to \bi{\Om}{r+1,s}, \ \
 d_T \colon \bi{\Om}{r,s} \to
\bi{\Om}{r,s+1}, \ \ \hbox{\rm and} \ \ \delta \colon \bi{\Om}{r,s} \to
\bi{\Om}{r-1,s+2} .
$$

% Since
%$d_{\mathcal F} ^ 2 = 0$ then
%$( \bi{\Om}{* ,0}, d_{\mathcal F})$ is a complex and its cohomology, called the foliated cohomology,
%we denote by $\hiru{\Om}{*}{M,\F}.$ It is the cohomology of
%the foliated manifold $(M,\F)$ considered as a topological
%space with two topologies.

\bigskip

 In this work, we shall use three types of
cohomologies.

\begin{itemize}
\item[(a)]
The {\em basic cohomology}
$\hiru{H}{*}{\mf}$ is the cohomology of the complex
$\hiru{\Om}{*}{\mf}$
of
{\em  basic forms}.
 A differential form $\om$
is basic when $i_{X}\om = i_{X}d\om =0$ for every vector field $X$
tangent to $\F$.
The complex  $\hiru{\Om}{*}{\mf}$ can be identified with the
complex of holonomy invariant forms on the transverse manifold $T$
- $\lau{\Om}{*}{\mathcal H}{T}.$
\end{itemize}

 \noindent In particular, if the foliation $\mathcal F$ is
developable and $ D \colon\tilde{M} \to  T $ its development with
connected fibres, $h \colon \pi_1(M) \to \Diff(T)$ - the
development representation, then the complex of basic forms
$\hiru{\Om}{*}{\mf}$ can be identified with the complex of
$h(\pi_1(M))$-invariant forms on $T.$

\begin{itemize}
\item[(b)] The {\em compactly supported  basic cohomology} $\lau{H}{*}{c}{\mf}$
 is the
 cohomology of the basic subcomplex  $\lau{\Om}{*}{c}{\mf}
 = \{ \om \in \hiru{\Om}{*}{\mf} \ |\  \hbox{ the support of
 $\om$ is compact}\}$.

\item[(c)] \label{zazpi} The {\em twisted basic cohomology} $\lau{H}{*}{\kappa}{\mf}$,
 relatively to the cycle $\kappa \in \hiru{\Om}{1}{\mf}$, is the
 cohomology of the basic complex $\hiru{\Om}{*}{\mf}$ relatively
 to the differential $ \ib{d}{\kappa}\om = d\om - \kappa \wedge \om$.
 This cohomology does not depend on the choice of the cycle: we
 have $\lau{H}{*}{\kappa}{\mf}\cong \lau{H}{*}{\kappa + df}{\mf}$
 through the isomorphism: $[\omega] \mapsto [e^f \omega]$. 

\end{itemize}

\prg{\bf Example} (E. Ghys, \cite{GH}).  Consider the unimodular
matrix $A=\left(%
\begin{array}{cc}
  1 & 0 \\
  1 & 1 \\
\end{array}%
\right)$  inducing a diffeomorphism of $\T^2 = \R^2/\Z^2 $. Let $\T^3_A$ be the
torus bundle over $\S^1$ determined by $A$ and $\mathcal F$ be the
flow obtained by suspending $A$. Then, the basic cohomology 
$\hiru{H}{2}{\T^3_A/ {\mathcal F}}$ is
infinite dimensional as basic forms
correspond to $A$-invariant forms on $\T^2$ - i.e.,
1-forms are of the form $f(x)dx,$ thus closed, and 2-forms are of
the form $u(x)dx \wedge dy.$

\medskip

 In \cite{RE_AJM}, Reinhart claimed that the basic cohomology of a Riemannian foliated closed manifold is
 finite dimensional and  satisfies the Poincaré duality property:

$$ \hiru{H}{k}{\mf}\cong \hiru{H}{n-k}{\mf}.$$

\noindent Soon it became apparent that the proof was not rigorous
and contained some gaps. For a long time it remained an open
problem. Only the beginnings of the eighties brought some striking
new developments. In addition to Sullivan's characterization of
taut foliations ({\it  A foliation is taut if and only if no
foliation cycle is the limit of boundaries of tangent chains.\/},
cf. \cite{SU}) and H. Rummler's thesis, \cite{RU}, A. Haefliger
published in 1980, perhaps the most influential result of this
theory, cf. \cite{HAE_JDG} - the property of being "taut" is a
transverse property, i.e. it depends only the properties of the
holonomy pseudogroup.  Then F. Kamber and Ph. Tondeur proved their
correct version of the Poincaré duality property for taut
Riemannian foliations on closed manifolds, cf. \cite{KT1,KT_AST}.
However, it was not until Y. Carrière's thesis (1981) that a
counterexample was found. The thesis presented the classification
of 1-dimensional tangentially orientable Riemannian foliations
(flows) on closed 3-manifolds. He found some flows which are not
defined by a Killing vector field, and these flows are
characterized by the property that the 2-dimensional basic
cohomology is trivial, cf. \cite{CA}.

\prg {\bf Flows}.
Let us begin with Carri\`ere's example.
Let $A$ be a matrix of $SL_2(\Z)$ with trace greater than 2 with two different eigenvalues
 $\lambda$ and $\frac{1}{\lambda}$ with two corresponding eigenvectors $v_1$ and $v_2,$ respectively.
 By $\T^3_A$ we denote the 3-manifold obtained by suspending $A,$ i.e. 
 it is a $\T^2$-fibre bundle over $\S^1.$
 In fact, it is obtained as the quotient space of $\T^2 \times \R$ by the equivalence relation
 generated by the identification of $(m,t)$ with $(A(m),t+1).$ The lines parallel to eigenvectors $v_1$
 and $v_2$ define $A$-invariant foliations (flows) $\Phi_1$ and $\Phi_2,$ respectively, on $\T^2.$ They,
 in turn, induce flows on $\T^3_A$ which we denote by the same letters. Each flow is dense in the tori
 which form the fibres of our 3-manifold $\T^3_A.$  One can show, cf. \cite{CA}, that
the flow $\Phi_2$ on $\T^3_A$ is a transversally Lie modelled
on the affine group $GA$ of the real line.

As the flow is transversally Lie, there exists a developing
mapping $D \colon {\mathbb R}^3 \rightarrow GA$ (${\mathbb R}^3$
is the universal covering of $\T^3_A$) such that the fibers of $D$
are the leaves of the lifted flow. Moreover, there exists a
homomorphism of groups $h\colon \pi_1(\T^3_A) \rightarrow GA$ and
its image is called the holonomy group $\Gamma$ of the foliation.
Global basic forms on $(\T^3_A, \Phi_2)$ correspond to
$\Gamma$-invariant forms on $GA,$ thus to $K =
\bar{\Gamma}$-invariant forms, where $K$ is the closure of the
group $\Gamma$ in $GA.$ Therefore the basic cohomology of the
foliated manifold $(\T^3_A, \Phi_2)$ is isomorphic the cohomology
of the complex of $K$-invariant forms on $GA.$ If we identify the
group $GA$ with the group ${\mathbb R}^2$ with the product given by
the formula $(t,s)(t',s') = (t+t', {\lambda}^ts' +s),$ then the
group $K$ can be identified with the group $\{ (n,s) \colon n \in
{\mathbb Z}, s \in {\mathbb R}\}.$

These consideration permit us to show that
$\hiru{H}{2}{\T^3_A/ \Phi_2} =0$.
To prove that fact we have to show that any $K$-left
invariant 2-form on $GA$ is exact. The 1-forms $\alpha = dt$ and
$\beta =\frac{ds}{{\lambda}^t}$ are left-invariant. A smooth
function if $K$ invariant if it does not depend on the variable
$s$ and $f(t) = f(t+1)$ for any real number $t.$ Hence a one form
$\omega$ is $K$-invariant iff $\omega = f\alpha +g\beta $ and both
functions $f$ and $g$ are $K$-invariant. A $K$-invariant 2-form
$\Omega$ can be written as $\Omega = h\alpha \wedge \beta$ where
$h$ is a $K$-invariant function. We have to demonstrate that for
any $K$-invariant function $h$ there exist $K$-invariant functions
$f$ and $g$ such that $d(f\alpha + g\beta ) = h \alpha \wedge
\beta , $ or equivalently $g'(t) + g(t)log\lambda =h(t)$ for any
real number $t.$

If we assume that $g(t) = {\lambda}^{-t}g_1(t), $ then we have to find
a function $g_1$ such that
$g'_1(t)\lambda^{-t}=h(t).$ But such a function is given by integration:

$$g_1(t) = c + \int_0^t {\lambda}^xh(x)dx$$

\noindent where $c$ is a real constant. And thus

$$g(t) = {\lambda}^{-t}( c + \int_0^t {\lambda}^xh(x)dx)$$

\noindent To get the invariance condition $g(t) = g(t+1)$ we need
$c= \frac{1}{\lambda -1}\int_0^1\lambda^xh(x)dx$ which is always
possible as $\lambda \neq 1.$ 

\medskip

Moreover, it is not difficult to show that the flow $\Phi_2$ is not isometric.

\medskip

This example should be seen in the light of the following
proposition, cf. \cite[ Proposition 6.6]{To}.

\bP
Let $\mathcal F$ be a flow defined by a nonsingular vector field
$V$ (with the normalized vector field $W= \frac{1}{\vert V \vert
} V$) on a Riemannian manifold $(M,\mu).$ Then the following
conditions are equivalent:
\begin{itemize}
\item[(i)] all leaves of $\mathcal F$ are minimal submanifolds of
$(M,\mu),$ i.e. the foliation is minimal;

\item[(ii)] the orbits of $V$ are geodesics;

\item[(iii)] $\theta (W) \ib{\chi}{\mu} = 0;$

\item[(iv)] $\nabla_WW= 0,$ where $\nabla$ is the Levi-Civita connection
of $(M,\mu).$
\end{itemize}
\eP

  The combined effort of H. Gluck, D. Sullivan, cf. \cite{GL,SU}, can be summarized
in the following theorem, cf. \cite[Proposition 6.7]{To}.
The equivalence of the fifth condition is due to Y. Carri\`ere,
cf. \cite{CA}.

\bT
 Let $\mathcal F$ be a flow given by the nonsingular vector field $V$ on an $m$-manifold $M$.
 Then the following conditions are equivalent:
\begin{itemize}
\item[(i)]   there exists a Riemannian metric on $M$ making the orbits of $V$ geodesics and $V$ of unit length;

\item[(ii)]  there exists a 1-form $\chi \in \hiru{\Om}{1}{M}$ such that $\chi (V) =1$ and $\theta (V) \chi = 0$;

\item[(iii)] there exists a 1-form $\chi \in \hiru{\Om}{1}{M}$ such that $\chi (V) =1$ and $i_Vd\chi =0$;

\item[(iv)]  there exists an (m-1)-plane subbundle $E \subset TM,$
complementary the flow such that $[V,X]$ is a section of $E$ for
any section $X$ of $E$,

 \item[(v)]  there exists a Riemannian metric on $M$ for which $V$ is a  Killing vector
 field.
 \end{itemize}
\eT

\prg{\bf Tautness and basic cohomology.}
In addition to all these  deep results a  simple remark seems to
point to a close relation between tautness and basic cohomology,
cf. \cite[Theorem 4.32]{To}.

\medskip

\bT
Let $\mathcal F$ be a transversally oriented Riemannian foliation
of codimension $n$ of a Riemannian manifold $(M,\mu)$ whose leaves are
minimal. Then the basic cohomology class of the transverse volume
$\nu$ form is non-trivial.
\eT

\pro Assume that there exists a basic form
$\alpha$ such that $d\alpha = \nu$. Let $\ib{\chi}{\mu}$ be the volume form
along the leaves of the foliation. Then

$$d ( \alpha \wedge \ib{\chi}{\mu} ) = d\alpha \wedge \ib{\chi}{\mu} + (-1)^{n-1} \alpha \wedge
d\ib{\chi}{\mu}$$

The minimality assumption implies that the form $d\ib{\chi}{\mu}$ is of
degree $(p-1,2)$, so as the form $\alpha$ is of degree $(0,n-1)$
the form $\alpha \wedge d\ib{\chi}{\mu}$ vanishes. Therefore the form $\nu
\wedge \ib{\chi}{\mu}$ which is a volume form of the manifold $M$ is exact,
a contradiction. \qed

\medskip

Moreover, in \cite{KT_AST}, F. Kamber and Ph. Tondeur proved that
the basic cohomology of a taut Riemannian foliation of a closed
manifold is finite dimensional and satisfies the PD property. The
above results and Haefliger's theorem, \cite{HAE_JDG}, which
assured that the existence of a Riemannian metric making all
leaves minimal is a transverse property made possible in 1982 the
formulation of the following conjecture by Y. Carrière first
expressed  for flows:

\medskip

\noindent{\bf Conjecture} { \it Let $\mathcal F$ be a Riemannian
foliation of a closed Riemannian manifold $(M,\mu).$ Then there
exists a (bundle-like) Riemannian metric making all leaves minimal
(i.e., the foliation is taut) iff the top dimensional basic
cohomology is non-trivial.}

\medskip

For flows the conjecture was solved by P. Molino and V. Sergiescu
in 1985, cf. \cite{MS}:

\bT
Let $\mathcal F$ be a Riemannian flow on a closed oriented
$m$-manifold $M$. Then there exists a Riemannian metric for which
$\mathcal F$ is an isometric flow iff the top dimensional basic
cohomology is non-trivial.
\eT

However, at that time the solution of the conjecture was far away.

\medskip

First, G. Hector and his students A. El Kacimi, and V. Sergiescu
proved that the basic cohomology of a Riemannian closed foliated
manifold $(M,g, {\mathcal F})$ is finite dimensional, cf.
\cite{EHS}, then they developed the Hodge theory, first studied by
F. Kamber and Ph. Tondeur in \cite{KT1,KT_AST}, for basic forms
and showed that the basic cohomology has the PD property iff the
top dimensional basic cohomology is non-trivial, cf. \cite{EH}.

\bT
Let $\mathcal F$ be a transversally oriented Riemannian foliation
on a closed oriented manifold $M.$  Then the following
two conditions are equivalent:
\begin{itemize}
\item[(i)] $\hiru{H}{n}{\mf} \neq 0;$

\item[(ii)] the basic cohomology $\hiru{H}{*}{\mf}$ satisfies the
Poincaré duality property.
\end{itemize}
\eT

This theorem together with Kamber-Tondeur's result mentioned at
the beginning of the subsection strongly hinted that Carrière's
intuition was correct. Finally in 1991, X. Masa, \cite{MA}, showed
the tautness is equivalent to the non-triviality of the
top-dimensional basic cohomology, solving  the
conjecture positively.

\bT
\label{Masa}
Let $\mathcal F$ be a transversally oriented
Riemannian foliation of a closed manifold $M$. Then there exists a
Riemannian metric on $M$ for which all leaves are minimal iff the
top-dimensional basic cohomology $\hiru{H}{n}{\mf}$ is
non-trivial.
\eT

% However, the theorem does not answers the question whether the
% metric can be chosen to be bundle-like.

To complete the story of basic cohomology let us mention that in
1993 A. El Kacimi and M. Nicolau proved  that the basic cohomology
of a closed Riemannian foliated manifold is a topological
invariant, cf. \cite{EN}.

\prg {\bf Mean curvature form}
In the story of the proof of the tautness conjecture a certain
1-form turned out to be of great importance.

\medskip

For a foliation $\mathcal F$ of a Riemannian manifold $(M,\mu),$ we
define the shape operator $W$ of the leaves using the natural
splitting of the tangent bundle

$$TM = T{\mathcal F} \oplus T{\mathcal F}^{\bot}.$$

\noindent In fact, for any section $Y$ of $T{\mathcal F}^{\bot}$
and any tangent vector field $X,$ we have

$$ W(Y)(X) = - \pi^{\bot}(\nabla_XY),$$
where $\pi^{\bot} \colon TM \to T{\mathcal F}^{\bot}$ is the orthogonal projection.

\noindent The trace of $W$ is linear in $Y,$ so it defines a section
of $T{{\mathcal F}^{\bot}}^{*}.$ We extend it to a global 1-form
$\ib{\kappa}{\mu}$ on $M:$

$$\ib{\kappa}{\mu} (X) =
\left\{
\begin{array}{ll}
    \trace W(X) & \hbox{if $X\in T\F^{\bot}$}\\[,3cm]
0& \hbox{ if
$X\in T{\mathcal F}$.}
\end{array}
\right.
$$

\noindent The 1-form is called the {\it mean curvature 1-form\/}
of $\mathcal F$ on the Riemannian manifold $(M,\mu).$

\medskip

If $f \colon (M_{1},\F_{1}) \to (M,\F)$ is a foliated
   imbedding
   between two Riemannian manifolds with $f(M_{1})$ saturated in $M$ and $\dim \F_1=\dim \F$, then 
   \begin{equation}
   \label{hamar}
   \hbox{$f^{*}\mu$ is a bundle-like
   metric on $(M_{1},\F_{1})$ and $f^{*}\ib{\kappa}{\mu} = \ib{\kappa}{f^{*}\mu}$.}
   \end{equation}
   So, 
   \begin{equation}
   \label{bost}
   \hbox{if $U$
 is an  open  subset of $M$ then
   $
   (\ib{\kappa}{\mu} )|_{U} = \ib{\kappa}{\mu|_{U}}.
   $
   }
   \end{equation}

%   There is a unique vector field $\ib{\tau}{\mu} $ with values in $T{\mathcal
%F}^{\bot} \subset TM$ such that

%$$\ib{\kappa}{\mu}  (V) =\mu(\ib{\tau}{\mu} , V)$$

%\noindent for any $V \in T{\mathcal F}^{\bot}.$

\medskip
This form is of particular interest. In \cite{KT1}, the authors
proved that if the form $\ib{\kappa}{\mu} $ is basic, then  it is closed. So
it defines a 1-basic cohomology class $[\ib{\kappa}{\mu} ]$ which proved to
be of importance in the study of taut foliations as, cf.
\cite{KT1},

\bP
\label{bat}
Let $\mathcal F$ be a Riemannian foliation
on a closed manifold $M$ with a bundle-like metric $\mu$ for which 
$\ib{\kappa}{\mu} $is basic and
$ [\ib{\kappa}{\mu} ]=0.$ Then bundle-like metric $\mu$ can be modified along
the leaves to a bundle-like metric $\mu'$ for which all leaves of
$\mathcal F$ are minimal.
\eP

\pro Since $ [\ib{\kappa}{\mu} ]=0,$ there exists a smooth basic
function $f$ on $(M/{\mathcal F})$ such that $ \ib{\kappa}{\mu} = df.$ Put
$\lambda = e^f$ and modify the metric $\mu$ as follows

$$ \mu' = {\lambda}^{\frac{2}{p}}\mu_{\mathcal F} \oplus \mu^{\bot}$$

\noindent where $p$ is the dimension of leaves, $\mu_{\mathcal F}$
and $\mu^{\bot}$ is the Riemannian metric induced on leaves of
$\mathcal F$ and the orthogonal subbundle, respectively. The
splitting is the splitting defined by the metric $\mu.$ The mean
curvature form $ \ib{\kappa}{\mu'} $ is equal to $ \ib{\kappa}{\mu} -
d \log\lambda =0.$ \qed

\medskip

Let $\mathcal F$ be a tangentially oriented foliation. We define
the characteristic form $\ib{\chi}{\mu},$  a $p$-form,
 as follows:

\medskip

\noindent for any p-tuple $(Y_1,...,Y_p)$ of vectors of $T_xM$

$$\ib{\chi}{\mu}(Y_1,...,Y_p) = \det (\mu(Y_i,E_j)_{ij})$$

\noindent where i,j = 1,...,p and $E_1,...,E_p$ is an oriented
orthonormal frame of $T_x{\mathcal F}.$

\medskip

There is a close relation between the characteristic form and the
mean curvature form. Namely, cf. \cite{RU},

\bT
Let $\mathcal F$ be a tangentially oriented foliation of a
Riemannian manifold $(M,\mu),$ $\ib{\chi}{\mu}$ its
characteristic form and $\ib{\kappa}{\mu}$ its mean curvature form. 
Then, for any
vector field $Y$ orthogonal to the foliation, we nave:

$$\theta (Y) \ib{\chi}{\mu} = - \ib{\kappa}{\mu} (Y) \ib{\chi}{\mu} +
\beta $$

\noindent where $\beta$ is a p-form of type $(p-1,1).$
\eT

As a corollary we get the following:

\medskip

\bC 
A tangentially oriented foliation $\mathcal F$ is taut iff
for any vector field $Y$ orthogonal to the foliation the form
$\theta (Y) \ib{\chi}{\mu}$ is of type $(p-1,1), $ which is
equivalent to the condition that $d\ib{\chi}{\mu}$ is of type
$(p-1,2), $ that is for any vector $Y$ and any vectors
$(Y_1,...,Y_p)$ tangent to the foliation

$$d\ib{\chi}{\mu}(Y,Y_1,...,Y_p) = 0.$$
\eC

The research into the tautness conjecture concentrated on the
study of the basic cohomology and the mean curvature form.

\medskip

The following theorem, cf. \cite{KT1,KT_AST}, gave further
evidence that the tautness, the mean curvature class  and the PD
property for basic cohomology are linked in some way. And the
result of A. El Kacimi and G. Hector suggested that the
non-vanishing of the top dimensional basic cohomology can be
related to the tautness of the foliation, i.e. the vanishing of
the mean curvature form.

\bT
Let $\mathcal F$ be a transversally oriented Riemannian foliation
foliation on a closed oriented manifold $M.$ Let $g$ be a
bundle-like metric with basic mean curvature form. Then the pairing
$\alpha \otimes \beta \rightarrow \int_M \alpha \wedge \beta
\wedge \ib{\chi}{\mu}$ induces a non-degenerate pairing

$$\hiru{H}{r}{\mf}\otimes  \lau{H}{n-r}{\ib{\kappa}{\mu}}{\mf}
\rightarrow {\mathbb R}$$

\noindent of finite-dimensional vector spaces.
\eT

\medskip

In the development of the theory Álvarez López's paper
\cite{Suso} of 1992 proved to be of great interest. In the
paper, Álvarez demonstrates that the space of smooth forms $\hiru{\Om}{}{M}$
on a foliated closed Riemannian manifold $(M,\mu,{\mathcal F})$ can
be decomposed as the direct sum of $\hiru{\Om}{}{\mf}$ of basic
forms and its orthogonal complement $\hiru{\Om}{}{\mf}^{\bot}.$
Therefore the mean curvature form $\ib{\kappa}{\mu}$ of $(M,\mu,{\mathcal F})$
can be decomposed into the basic component $\ib{\kappa}{\mu,b}$ and the
orthogonal one. The 1-form $\ib{\kappa}{\mu,b}$ is closed and it defines the
1-basic cohomology class $\kappab = [\ib{\kappa}{\mu,b}]$, which does not depend on $\mu$. \label{bi} Moreover, Álvarez proves that
any form cohomologous to $\ib{\kappa}{\mu,b}$ (in the complex of basic
forms) can be realized as the basic component of the mean
curvature form of some bundle-like metric of $\mathcal F$ with the
same transverse Riemannian metric. Additionally, one can verify
that changing the orthogonal  complement of $\mathcal F$ does not
change the form $\ib{\kappa}{\mu,b}.$

As an application, the assumption of the orientability of $M$ in the 
original formulation of Theorem   \ref{Masa} is removed.

For some time the condition that the mean curvature form is basic
seemed to be a major obstacle to the existence of such a
Riemannian metric. But at last in 1995 D. Domínguez published
his theorem stating that, \cite{DO1,D},

\bT
\label{TD}
 Let $\mathcal F$ be a Riemannian foliation on a
closed manifold $M$. Then there exists a bundle-like metric for
$\mathcal F$ for which the mean curvature form is basic. 
\eT
These Riemannian metrics are very important in the remaining part of the
paper. Therefore, a  bundle-like metric for which the mean curvature form is basic
we call a {\em $D$-metric}.

 This Theorem together with Proposition \ref{bat} ensures that the ``taut"
 Riemannian  metric can be chosen to be a $D$-metric. In the sequel we shall use the following fact:  
 \begin{equation}
 \label{rest}
 \hbox{if $U$ is a
saturated open subset of $M$ such that 
 $\mu|_{U}$ is a D-metric then $\ib{\kappa}{\mu,b}|_{U} =
 \ib{\kappa}{\mu|_{U}}$.
 }
 \end{equation}

The final characterization of taut Riemannian foliations of closed
manifolds can be summarized in  the following theorem, cf.
\cite[7.56]{To}:

\bT Let $\mathcal F$ be a transversally oriented
Riemannian foliation foliation on a closed oriented Riemannian
manifold $(M,\mu).$  Then
 $\lau{H}{n}{\ib{\kappa}{\mu}} {\mf}\cong {\mathbb R}$ and the
following conditions are equivalent:

\begin{itemize}
\item[(i)] $\hiru{H}{n}{\mf} \cong {\mathbb R},$

\item[(ii)] ${\mathcal F}$ is taut;

\item[(iii)] $\kappab= 0;$

\item[(iv)]  $\lau{H}{0}{\ib{\kappa}{\mu}}{\mf}= \R.$ 
\end{itemize}
Moreover, then the basic cohomology of the foliated manifold $(M/{\mathcal F})$ has the Poincaré duality property.
\eT

\prg{\bf Open manifolds}
The theory has not been well-developed for open manifolds. We have
a fine and very general version of Poincaré duality theorem
published by V. Sergiescu in 1985, cf. \cite{Se}. Then in 1997,
Cairns and Escobales presented a very interesting example, cf.
\cite{CE}, of a Riemannian foliation on an open manifold for which
the mean curvature form is basic but not closed.

\prgg{\bf The SRFs }. A singular Riemannian
foliation\protect\footnote{For the notions related
with singular Riemannian foliations we refer the reader to
\cite{BM,Mo0,Mo,Mo1}.}
(SRF for short) on a connected manifold $X$ is
a partition $\mathcal{K}$ by connected immersed
submanifolds, called {\em leaves}, verifying the following properties:

\begin{itemize}
\item[I-] The module of smooth vector fields tangent to the leaves is
transitive on each leaf.

\item[II-] There exists a Riemannian metric $\nu$ on $N$, called {\em
adapted metric},  such that each geodesic
that is perpendicular at one point to a leaf remains perpendicular to every
leaf it meets.
\end{itemize}
The first condition implies that $(X,\mathcal{K})$ is a singular
foliation in the sense of \cite{St} and \cite{Su}. Notice that the
restriction of $\mathcal{K}$ to a saturated open subset produces a
SRF. Each (regular) Riemannian foliation (RF in short) is a SRF, but
the first interesting examples are the following:
\begin{itemize}
\item[-] The orbits of the action by isometries of a Lie group.
\item[-] The closures of the leaves of a regular Riemannian foliation.
%     \item[-] The sheets of a totally geodesic foliation.
\end{itemize}

\prgg{\bf Stratification}.  Classifying the points of
  $X$ by the dimension of the leaves one gets a stratification
  $\SK$ of
  $X$ whose elements are called  strata. The restriction of
  $\mathcal{K}$ to a stratum $S$ is the RF
  $\mathcal{K}_{S}$. The strata are ordered by: $S_{1} \preceq S_{2}
  \Leftrightarrow S_{1}\subset \overline{S_{2}}$. The minimal (resp.
  maximal) strata are the closed strata (resp. open strata). We shall
  denote by  $S_{_{min}}$ the union of the closed strata.
  Since $X$ is
  connected, there is just one open stratum, denoted
  $R_{\mathcal{K}}$. It is a dense subset.
  This is the {\em regular stratum}, the other strata are the {\em
  singular strata}.

  The  depth  of $\SK$, written $ \depth  \SK$, is defined to be the largest
 $i$ for which there
 exists a chain of strata $S_0  \prec S_1 \prec \cdots \prec S_i$. So,
 $ \depth  \SK= 0$ if and only if the foliation
 ${\mathcal  K}$ is regular. The  {\em depth}  of  a
 stratum $S \in \SH$ , written $ \depth_{\HH} S$, is defined to be the largest
  $i$ for which there
  exists a chain of strata $S_0  \prec S_1 \prec \cdots \prec S_i= S$.
\medskip

The basic cohomology of such foliations on closed manifolds is
finite dimensional and it is a topological invariant, cf.
\cite{W}. However, as far as the tautness property is
concerned the situation is totally different.

\prgg {\bf Example.}
Let us consider the isometric action
      $
      \Phi \colon {\mathbb R} \times \S^{2d+2} \to \S^{2d+2}
      $
      given by the formula
      $$
      \Phi (t, (z_{0}, \ldots , z_{d},x)) = (e^{a_{0} \pi it} \cdot z_{0}, \ldots ,
      e^{a_{d}\pi it}\cdot z_{d} ,x),
      $$
      with $(a_{0}, \ldots , a_{d}) \ne (0,\ldots , 0)$. Here,
      $
      \S^{2d+2} = \{ (z_{0}, \ldots , z_{d},x) \in \mathbb{C}^{d+1} \times {\mathbb R} \  |  \ |z_{0}|^{2} + \cdots +
      |z_{d}|^{2}+ x^{2}= 1 \}$.
      There are two singular strata: the north pole $S_{1}= (0, \ldots ,
      0,1)$ and the south pole $S_{2} = (0,\ldots,-1)$.
      The regular stratum is $\S^{2d+1} \times ]-1,1[$. Let  $r$ be
      the variable of $]-1,1[$.
      The basic cohomology $\hiru{H}{*}{\S^{2d+2}/\mathcal F}$ of
      the foliation is

\begin{tiny}
      \begin{center}
      \begin{tabular}{|c|c|c|c|c|c|c|c|c|c|}\hline
          &&&&&&&&\\
       $i=0$ & $i=1$ & $i=2$  & $i=3$ & $i=4$ & $i=5$  & $i= \cdots $ &
       $i=2d$ & $i=2d+1$    \\[,1cm] \hline
       &&&&&&&&\\
        $1$&  0 & 0  & $ [dr
      \wedge e] $
      &0 & $[dr \wedge e^2] $ & $\cdots$ & 0 &  $ [dr \wedge e^d] $    \\[,1cm]

\hline
      \end{tabular}
      \end{center}
\end{tiny}
      \medskip

      \noindent where $e \in \lau{\Omega}{2}{\bar{2}}{\S^{2d+2}/ \mathcal F}$
      is an  Euler
      form.
% (cf. \cite{HS}).

The top dimensional basic cohomology group is
      isomorphic to $\mathbb R $, but this cohomology does not have the Poincaré
      duality property  in spite of the fact that the flow is
      isometric. And, of course, the foliation is not minimal for
      any adapted (bundle-like) Riemannian metric.

Moreover, in \cite{MW}, the authors proved that a singular
 foliation on a closed manifold admitting an adapted Riemannian
metric for which all leaves are minimal must be regular. These
fact have led us to study closer singular Riemannian foliations.
We have introduced basic intersection cohomology  in view to
recover some kind of Poincaré duality, cf. \cite{SW1,SW2,PSW1}.
We hope that soon we will complete our task and demonstrate the
perverse version of the Poincaré duality property for basic
intersection cohomology for singular Riemannian foliations of
closed manifolds. In his thesis, \cite{JI}, written under the supervision of
M. Saralegi and M. Macho, J.I. Royo Prieto demonstrated,  among other
results on singular Riemannian flows,
the Poincaré duality for basic intersection cohomology and the singular version of
the Molino-Sergiescu theorem. Inspired by these results,  we have started to investigate the possible
generalizations to the SRF case and at the same time we
have found that our research gives some interesting insights into the
problem on non-compact manifolds, cf. \cite{PSW1,RSW}.  The second part
of this work is concerned with this problem.

We complete the section with the presentation of the BIC for the
above example, in which the PD property can be easily seen.

If we consider the BIC of our  example the picture changes.
      The following table presents the BIC
      $\lau{IH}{*}{\bar{p}}{\S^{k =2d+2}/\mathcal{F}}$ for the constant
      perversities:

\vspace{1cm}

{\begin{tiny}
 \begin{center}
\kern-1truecm\begin{tabular}{|c|c|c|c|c|c|c|c|c|c|c|c|c|}\hline
          &&&&&&&&&&&\\
     $ i \! = \!$ & $ 0$ & $ 1$ & $ 2$  & $ 3$ & $ 4$ & $ 5$  &$6$& $7$& $  \cdots $ & $
      k \! - \! 2$ & $    k \! - \!  1$
      \\[,1cm] \hline
      &&&&&&&&&&&\\
       $\bar{p} <\bar{0}$&0&$[dr]$ & 0&$[e \wedge dr]$ &  0 & $[e^2 \wedge dr]$
       &  0 & $[e^3 \wedge dr]$ &$\cdots$ &0& $[e^d \wedge dr]$
       \\[,1cm]  \hline
       &&&&&&&&&&&\\
       $\bar{p}  \! = \!\bar{0},\bar{1}$ &  $1$ & 0 &  0 &  $[e \wedge dr]$ &  0 &$
           [e^2 \wedge dr]$&  0 & $[e^3 \wedge dr]$ & $\cdots$  & 0&$[e^d \wedge dr]$
           \\[,1cm]  \hline
           &&&&&&&&&&&\\
      $\bar{p}  \! = \!\bar{2},\bar{3}$&0&0 & [e]& 0 &  0 & $[e^2 \wedge dr]$&  0 & $[e^3 \wedge dr]$ &$\cdots$
      &0&$[e^d \wedge dr]$
      \\[,1cm]  \hline
      &&&&&&&&&&&\\
      $\bar{p}  \! = \!\bar{4},\bar{5}$ &  $1$ & 0 &  $[e]$ & $0$ &  $[e^{2}]$ & $0$&  0 & $[e^3 \wedge dr]$
      &$\cdots$ &0&$[e^d \wedge dr]$
      \\[,1cm]  \hline      &&&&&&&&&&&\\
      $\cdots$ & $\cdots$ & $\cdots$ & $\cdots $& $\cdots$ & $\cdots$& $\cdots$ & $\cdots$ &
      $\cdots$  & $\cdots$& $\cdots$& $\cdots$\\[,1cm]
      \hline       &&&&&&&&&&&\\
      $\bar{p}  \! = \!\bar{k  \! \! - \!   \!  4}, \bar{k \! \! - \! \!
      3}$ &  $1$ & 0 &  $[e]$ &
      0 & $[e^{2}]$& 0 & $[e^{3}]$&$0$ &$\cdots$ &0&$[e^d \wedge dr]$
      \\[,1cm] \hline       &&&&&&&&&&&\\
      $\bar{p} \geq\bar{k \! - \! 2} $ &  $1$ & 0 &  $[e]$ &
      0 & $[e^2]$ &  0 & $[e^{3}]$&0& $\cdots$&$[e^{d}]$& 0\\[,1cm]  \hline
      \end{tabular}
    \end{center}
\end{tiny}
}

\bigskip

      We notice that the top dimensional basic cohomology group is
  isomorphic either to 0 or ${\mathbb R} $. These cohomology groups are finite
  dimensional. We recover the Poincaré  duality in the perverse sense:

      $$
      \lau{IH}{*}{\bar{p}}{\S^{k}/\mathcal{F}} \cong \lau{IH}{k-1-*}{\bar{q}}{\S^{k}/ \mathcal{F}}
      $$
for two complementary perversities: $\bar{p} + \bar{q} =  \bar{t} =\overline{k-3}$.

\medskip

\bT
 Let $M$ be a connected closed manifold endowed with
an SRF $\mathcal F$.
 If $ \ell = codim_{M} {\mathcal F}$ and $\bar{p}$ a perversity on $(M/{\mathcal F})$,
 then
 $$
 \lau{IH}{\ell}{\bar{p}}{\mf} = 0 \hbox{ or } {\mathbb R}.
 $$
\eT

\medskip

\bC
 Let $M$ be a connected compact
manifold endowed with an SRF $\mathcal F$. Let  us suppose that
$\mathcal F$ is transversally orientable. Consider  $\bar{p}$ a
perversity on $(M, {\mathcal F})$ with   $\bar{p}\leq \bar{t}$. If
$\ell= codim_{M} \mathcal F$, then the two following statements are
equivalent:

\begin{itemize}
\item[(1)]  The foliation ${\mathcal F}_{R}$ is taut, where $R$ is the
regular stratum of $(M,{\mathcal F})$;

\item[(2)] The cohomology group $ \lau{IH}{\ell}{\bar{p}}{\mf}$ is
${\mathbb R} $.
\end{itemize}
\noindent Otherwise, $ \lau{IH}{\ell}{\bar{p}}{\mf}= 0$.
\eC

\medskip

The BIC of a conical foliation  $\mathcal F$ defined on $M$ by an
isometric  action of an abelian Lie group on an oriented manifold $M$
verifies the Poincaré duality:
\begin{equation}
 \lau{IH}{\ell}{\bar{p}}{\mf} \cong
\lau{IH}{\ell -*}{\bar{q},c}{\mf}.
\end{equation}
Here  $\ell= codim_{M} \mathcal F$ and the  two  perversities $\bar{p}$ and $\bar{q}$
are complementary.

\medskip

\noindent {\bf Note:} Due to the limited space we could dedicate to
this overview of the problem we have not mentioned many partial
results, (e.g. \cite{Suso,KT,H2}), and some reviews papers (e.g.
\cite{C1,SER_A}).

\section{Geometrical preliminaries.}
We present in this section the  kind of foliations we are going to
use in this work: the CERFs.
A CERF is essentially a Riemannian foliation defined on a non-compact
manifold which is imbeddable in a closed manifold in a nice way.

\prg {\bf The CERFs.}  We shall consider in this work a particular
case of Riemannian foliations defined on non-compact manifolds. They
have an outside compact manifold (zipper) and an inside compact submanifold (reppiz).
Consider a manifold $M$ endowed with a
Riemannian foliation $\F$.

\smallskip

A {\em zipper} of $\F$ is a closed manifold $N$ endowed with a (regular) Riemannian
foliation $\HH$ verifying the following properties:

\begin{itemize}
    \item[(a)] The manifold $M$ is a saturated open subset of $N$ and
    $\HH_{M} = \F$.
    \end{itemize}
    The open subset $M$ is also $\overline{\F}$-saturated. Thus,
    the closure $\overline{L}$ of a leaf $L \in \F$ is compact.
\smallskip

A {\em reppiz} of $\F$  is a saturated open subset $U$ of $M$ verifying the following
    properties:

    \begin{itemize}
    \item[(b)] the closure $\overline{U}$ (in $M$) is compact.
    \item[(c)] the inclusion $U \hookrightarrow M$ induces the isomorphism
$\hiru{H}{*}{U/\F} \cong \hiru{H}{*}{\mf}$.
    \end{itemize}
      It is not true that any saturated open subset of $M$ is a reppiz.
Just consider $M= \sbat$ endowed with  the pointwise
    foliation  and take $U = \sbat\backslash \{ (\cos (2\pi/n),\sin
    (2\pi/n)) \ / \ n \in \N\backslash \{ 0 \} \}$.

    \smallskip

    We say that $\F$
is a {\em {\rm C}ompactly
{\rm E}mbeddable {\rm R}iemannian {\rm F}oliation } (or
CERF)\footnote{The definition of CERF given in \cite{RSW} is more
restrictive: see Proposition \ref{SC}.} if
$(M,\F)$ possesses a zipper and a reppiz.
When $M$ is closed, then $(M,\F)$ is clearly a CERF, being $M$
itself a zipper and a reppiz.
Neither the zipper nor the reppiz are  unique.

The main example of a CERF is given by the strata of a singular
Riemannian foliation defined on a closed manifold. This family will
be treated in the next Section. The interior of a Riemannian foliation
defined on a
manifold with boundary is a CERF when the foliation is tangent to the
boundary; we can consider the double of the manifold as a zipper.
When the foliation is transverse to the boundary then the foliation is
not a CERF.

\bigskip

We present now some geometrical tools we shall use for the study 
of a SRF $(X,\mathcal{K})$.

  \prg {\bf Tubular neighborhood.}
  A singular stratum $S \in \SK$ is a proper submanifold of the
  Riemannian manifold $(X,\nu)$.
  So, it possesses a tubular neighborhood $(T_S,\tau_S,S)$.
  Recall that associated with this neighborhood there are the following smooth maps:
  \begin{itemize}
      \item[+] The {\em radius map} $\rho_S \colon T_S \to [0,1[$ defined
      fiberwise by
  $z\mapsto |z|$. Each $t\not= 0$ is a regular value of the $\rho_S$.
  The pre-image $\rho_S^{-1}(0)$ is $S$.
       \item[+] The {\em contraction} $H_S \colon T_S \times [0,1] \to T_S$
       defined fiberwise by  $(z,r ) \mapsto  r \cdot z$. The
       restriction
       $(H_S)_t \colon T_S \to T_S$ is an
       imbedding for each
       $t\not= 0$   and $(H_{S})_0 \equiv \tau_S$.
  \end{itemize}

  \nt These maps verify $\rho_S(r \cdot z) = r  \rho_S(z)$.
This tubular neighborhood can be chosen verifying the two following important
properties (cf. \cite{Mo}):

  \Zati Each  $(\rho_S^{-1}(t),\mathcal{K})$ is a
      SRF, and

       \zati Each  $(H_S)_{t} \colon  (T_S,\mathcal{F}) \to (T_S,\mathcal{F}) $ is a
       foliated map.

       \medskip

       \nt We shall say that $(T_S,\tau_S,S)$   is a {\em foliated tubular
neighborhood} of $S$.

The hypersurface $D_S = \rho_S^{-1}(1/2)$ is the {\em
     core} of the tubular neighborhood. We have  the equality
     $\depth \SKDS= \depth
     \SKTS -1$. Notice that the map
     \begin{equation}
         \label{difeo}
         \mathfrak{L}_{S}\colon (D_{S}\times
     ]0,1[,\K \times \I) \to ((T_{S}\backslash S),\K),
    \end{equation}
    defined by
     $\mathfrak{L}_{S}(z,t) = H_S(z,2t)$, is a foliated diffeomorphism.

\medskip

A family of foliated tubular neighborhoods $\{T_{S } \ | \ S \in
\SF^{^{sin}}\}$ is a {\em foliated Thom-Mather system} of $(N,\HH)$
if the following conditions are verified.

\begin{itemize}
    \item[]
    \begin{itemize}
    \item[(TM1)] For each pair of singular strata $S, S'$ we have
$$
T_{S} \cap T_{S'} \ne \emptyset  \Longleftrightarrow  S \preceq S'
\hbox{ or } S' \preceq S.
$$
\end{itemize}
\end{itemize}
\nt Let us suppose that $S'\prec S$. The two other conditions are:

\begin{itemize}
    \item[]
    \begin{itemize}

\item[(TM2)] $T_{S} \cap T_{S'} = \tau_{S}^{-1}(T_{S'} \cap S)$.

\item[(TM3)] $
   \rho_{S'} = \rho_{S'} \rondp \tau_{S}
   \hbox{ on }  T_{S} \cap T_{S'} .$
\end{itemize}
\end{itemize}

    We have seen in \cite{RSW} that
   each closed manifold endowed with a SRF possesses a foliated
    Thom-Mather system.  We fix for the sequel of this work a such foliated Thom-Mather system.

     \prg {\bf Blow up}.  Molino's blow up  of a SRF produces a
new SRF of the same generic dimension but with smaller
depth (see \cite{Mo} and also \cite{SW2},\cite{RSW}).
The main idea is to replace each point of the closed strata by its
link (a sphere).

\smallskip

In fact, given a SRF $(X,\mathcal{K})$ with $\depth \SK> 0$, there
exists another SRF $(\wh{X},\wh{\mathcal{K}})$ and a continuous map
$\mathfrak{L} \colon \wh{X} \to X$, called {\em blow up} of
$(X,\mathcal{K})$, verifying:
\begin{itemize}
    \item[-] $\depth \SKh = \depth \SK -1$.
    \item[-] there exists a commutative diagram
    $$
    \begin{picture}(150,70)(00,0)

    \put(0,60){\makebox(0,0){$ \mathfrak{L}^{-1} (X\backslash
    S_{_{min}}) $}}
    \put(150,60){\makebox(0,0){$ (X\backslash
    S_{_{min}}) \times \{ -1,1\}$}}
    \put(80,00){\makebox(0,0){$ X\backslash
    S_{_{min}}$}}

    \put(10,50){\vector(1,-1){40}}
    \put(145,50){\vector(-1,-1){40}}
    \put(40,60){\vector(1,0){56}}

    \put(68,69){\makebox(0,0){$f_{0}$}}
    \put(15,30){\makebox(0,0){$\mathfrak{L}$}}
    \put(160,30){\makebox(0,0){\small{projection}}}
    \end{picture}
    $$
    where $f_{0} \colon (\mathfrak{L}^{-1} (X\backslash
    S_{_{min}}) ,\wh{\mathcal{K}}) \to (X\backslash
    S_{_{min}} \times \{-1,1\} , \mathcal{K} \times \mathcal{I}) $ is
    a foliated diffeomorphism. Here, $\mathcal{I}$ denotes the
    foliation by points.

    \item[-] for each minimal (closed) stratum $S_{_{c}}$, there exists a commutative diagram
    $$
    \begin{picture}(150,70)(00,0)

        \put(10,60){\makebox(0,0){$ \mathfrak{L}^{-1} (T_{S_{_{c}}})$}}
    \put(140,60){\makebox(0,0){$ D_{S_{_{c}}}\times ]-1,1[$}}
    \put(80,00){\makebox(0,0){$ T_{S_{_{c}}}$}}

    \put(20,50){\vector(1,-1){40}}
    \put(135,50){\vector(-1,-1){40}}
    \put(40,60){\vector(1,0){56}}

    \put(68,69){\makebox(0,0){$f_{S_{_{c}}}$}}
    \put(25,30){\makebox(0,0){$\mathfrak{L}$}}
    \put(140,30){\makebox(0,0){$\mathfrak{L}_{S_{_{c}}}$}}
    \end{picture}
    $$
    where $f_{S_{_{c}}} \colon (\mathfrak{L}^{-1} (T_{S_{_{c}}}) ,\wh{\mathcal{K}})
    \to (D_{S_{_{c}}}\times ]-1,1[ , \mathcal{K} \times \mathcal{I}) $ is
    a foliated diffeomorphism and the map $\mathfrak{L}_{S_{_{c}}}$ is defined
    by $\mathfrak{L}_{S_{_{c}}}(z,t) = H_{S_c}(z,2 |t|)$.
    Notice that
    $f_{S_{_{c}}} \colon (\mathfrak{L}^{-1} (S_{_{c}}) ,\wh{\mathcal{K}})
        \to (D_{S_{_{c}}}\times \{ 0 \} , \mathcal{K} \times \mathcal{I}) $
        is also a foliated diffeomorphism.
    \end{itemize}

\medskip

The stratification induced by $\wh{\mathcal{K}}$ can be described as
follows. For each non minimal stratum $S \in \SK$ there exists a unique
stratum
$S^{^{\mathfrak{L}} }\in \SKh$ with $\mathfrak{L}^{-1}(S) \subset S^{^{\mathfrak{L}} } $,
and we have
$$
  \SKh = \{ S^{^{\mathfrak{L}} } \ / \ S \in \SK \hbox{ and } S \cap
S_{_{min}} = \emptyset\}.
$$
In fact,
$$
\left\{ \begin{array}{rcl}
f_{0}\left( S^{^{\mathfrak{L}} } \cap \mathfrak{L}^{-1}(X\backslash
S_{_{min}})\right) &= & S \times \{-1,1\} \hbox{ and } \\[,2cm]
f_{S_{_{c}}}\left( S^{^{\mathfrak{L}} } \cap
\mathfrak{L}^{-1}(T_{S_{_{c}}} )\right) &=  &
(S \cap D_{S_{_{c}}}) \times ]-1,1[
\hbox{ if $S_{_{c}}$ is a closed stratum with $S_{_{c}} \preceq S$}.
\end{array}
\right.
$$

\medskip

The CERFs and the SRFs are related by the following result.

\bp
\label{SC}
Let $X$ be a  closed manifold endowed with a SRF $\mathcal{K}$. For
any
stratum $S$ of $\SK$ the foliation $\mathcal{K}_{S}$ is a CERF.
\ep
\pro
When $S$ is a closed stratum it suffices to take the zipper
$(S,\mathcal{K})$  and the reppiz $S$. Consider now the case
where $S$ is not closed (minimal). We proceed in two steps.

\medskip

\underline{\em A zipper for $(S,\mathcal{K})$}.
Proceeding by induction on $\depth \SK$ we know that there
exists a zipper $(N,\HH)$ of $(S^{^{\mathfrak{L}} }
,\wh{\mathcal{K}})$. Since the map $\xi \colon (S,\mathcal{K}) \to (S^{^{\mathfrak{L}} }
,\wh{\mathcal{K}})$, defined by $x \mapsto f_{0}^{-1}(x,1)$,  is an
open  foliated imbedding we can identify
$(S,\mathcal{K}) $ with its (open) image $(\xi(S)
,\wh{\mathcal{K}})$. So, the foliated manifold
$(N,\HH)$ is a zipper of $(S,\mathcal{K})$.

\smallskip

\underline{\em A reppiz for $(S,\mathcal{K})$.}
For each $i\in \{0, \ldots , s -1\}$, where $s = \depth _{\HH} S$,  we
denote by :
\begin{itemize}
    \item[-] $\Sigma_{i} = \cup \{ S' \in \SH \ | \ \depth _{\HH} S' \leq
 i\}$,
\item[-] $T_{i}$ the union of the disjoint tubular neighborhoods  $\{ \ib{T}{S'} \ / \    \ib{T}{S'} \subset \Sigma_{i}\menos \Sigma_{i-1}\}$,
\item[-] $\rho_{i} \colon T_{i}\to [0,1[$ its radius function, and
    \item[-] $D_{i} = \rho_i^{-1}(0)$ the core of $T_{i}$.

\end{itemize}

The family
$\{S\cap T_{0},S\backslash \rho_{0}^{-1}([0,7/8])\}$ is a saturated
open covering of $S$.
The inclusion
$$I \colon ((S \cap T_{0} )\backslash
\rho_{0}^{-1}([0,7/8]),\mathcal{K})
\hookrightarrow (S\cap T_{0},\mathcal{K}) $$ induces an
isomorphism for the basic cohomology. This comes from the fact that
the inclusion $I$  is foliated
diffeomorphic to the inclusion
$$
J \colon ((S \cap D_{0})\times ]7/8,1[ , \mathcal{K}\times \I) \hookrightarrow
((S \cap D_{0} )\times ]0,1[ , \mathcal{K} \times \I)
$$
(cf. \refp{difeo} and $S \cap \Sigma_{0}= \emptyset$).
From the Mayer-Vietoris sequence (see for example \cite{RSW}) 
we conclude that the inclusion
$
S\backslash \rho_{0}^{-1}([0,7/8])
\hookrightarrow
S
$
induces the isomorphism
$$
\hiru{H}{*}{S/\mathcal{K}} \cong
\hiru{H}{*}{\left(S\backslash \rho_{0}^{-1}([0,7/8])\right)/\mathcal{K}}.
$$

The family
$$
\{ T_{S'} \menos  \rho_{0}^{-1}([0,7/8]) \ | \ S' \in \SF  ,
\depth_{\HH} S'> 0 \}
$$
is a foliated Thom-Mather system of
$(S\menos\rho_{0}^{-1}([0,7/8]),\HH)$ (cf.
\cite[(1.6)]{RSW}).
The same
previous argument applied to the stratum $S\backslash \rho_{0}^{-1}([0,7/8])$
gives
$$
\hiru{H}{*}{\left( S\backslash \rho_{0}^{-1}([0,7/8]) \right)/\mathcal{K}} \cong
\hiru{H}{*}{\left( \left( S\backslash \rho_{0}^{-1}([0,7/8]) \right) \backslash \rho_{1}^{-1}([0,7/8]) \right)/\mathcal{K}}.
$$
This procedure leads us to
\begin{eqnarray*}
\hiru{H}{*}{S/\mathcal{K}}
&\cong&
\hiru{H}{*}{(S\backslash \rho_{0}^{-1}([0,7/8]) )/\mathcal{K}}
\cong
\hiru{H}{*}{(S\backslash (\rho_{0}^{-1}([0,7/8]) \cup
\rho_{1}^{-1}([0,7/8]) ))/\mathcal{K}}
\cong
\cdots \\
&\cong&
\hiru{H}{*}{(S\backslash (\rho_{0}^{-1}([0,7/8]) \cup  \cdots \cup
\rho_{s-1}^{-1}([0,7/8]) ))/\mathcal{K}}.
\end{eqnarray*}
Take $$
U = S\backslash (\rho_{0}^{-1}([0,7/8]) \cup  \cdots \cup
\rho_{s-1}^{-1}([0,7/8]) ),
$$
which is an open saturated subset included
on $S$. By construction, the inclusion $U \hookrightarrow S$ induces
the isomorphism
$
\hiru{H}{*}{S/\mathcal{K}}
\cong
\hiru{H}{*}{U/\mathcal{K}}.
$
This gives (a).

Consider $K = S\backslash (\rho_{0}^{-1}([0,1/8[) \cup  \cdots \cup
\rho_{s-1}^{-1}([0,1/8[) )$, which is a subset of $S$ containing $U$.
We compute its closure in $S$:
\begin{eqnarray*}
    \overline{K} &=& \overline{S\backslash (\rho_{0}^{-1}([0,1/8]) \cup  \cdots \cup
\rho_{s-1}^{-1}([0,1/8]) )} \subset
\overline{S}\backslash \left((\rho_{0}^{-1}([0,1/8[))^{\hbox{°}} \cup  \cdots \cup
(\rho_{s-1}^{-1}([0,1/8[))^{\hbox{°}}\right)\\
&=& \overline{S}\backslash \left(\rho_{0}^{-1}([0,1/8[) \cup  \cdots \cup
\rho_{s-1}^{-1}([0,1/8[)\right) = S\backslash (\rho_{0}^{-1}([0,1/8[) \cup  \cdots \cup
\rho_{s-1}^{-1}([0,1/8[) ),
\end{eqnarray*}
since $\overline{S}\backslash S \subset \Sigma_{s-1} =
\rho_{0}^{-1}(\{0 \}) \cup  \cdots \cup
\rho_{s-1}^{-1}(\{ 0 \}) $.
This implies that $K$ is a closed subset of $S$ and therefore
compact.
This gives (b).
\qed

\prg {\bf Basic cohomology}.
As in the regular case, the {\em basic cohomology}
$\hiru{H}{*}{X/\mathcal{K}}$ is the cohomology of the complex
$\hiru{\Om}{*}{X/\mathcal{K}}$
of
{\em  basic forms} (cf. \cite{W}).
 A differential form $\om$
is basic when $i_{X}\om = i_{X}d\om =0$ for every vector field $X$
tangent to $\F$.

Associated to a covering $\{ U , V \}$  of $X$ by saturated open
 subsets we have the Mayer-Vietoris short exact sequence
 $$
 0 \to (\hiru{\Om}{*}{X/\mathcal{K}} ,d)\to
 (\hiru{\Om}{*}{U/\mathcal{K}},d) \oplus
 (\hiru{\Om}{*}{V/\mathcal{K}},d)  \to
 (\hiru{\Om}{*}{(U \cap V)/\mathcal{K}},d)  \to 0,
 $$
 where the maps are defined by restriction
 (the same proof of \cite{RSW} for the regular case works).

\section{Tautness in the non-compact case.}

We prove in this section that the previous cohomological characterizations of
the tautness of a RF $\F$ are still valid when
the manifold is non-compact but
the foliation $\F$ is a CERF.

\medskip

We fix for the rest of this section a CERF $\F$ defined on a
manifold $M$. We also fix a zipper $(N,\HH)$ and a reppiz $U$.

\prg{\bf Tautness class of \hbox{\boldmath {$\F$}}.}
Since $N$ is compact we get from Theorem \ref{TD}
that $M$ possesses a D-metric $\mu$. The {\em tautness class} of $(M,\F)$ is
the class $\kappab = [\ib{\kappa}{\mu} ] \in \hiru{H}{1}{\mf}$ (cf. page
\pageref{bi}).
This class is well defined since:

\bP
\label{2D}
Two D-metrics  on $(M,\F)$ define the same tautness class.
\eP
\pro
Fix a zipper $(N,\HH)$ and a reppiz $U$.
Since $N$ is compact then the tautness class $\kappab_{N}$ is well
defined.
Let $\mu$ be a D-metric on $M$. The key point of the proof is to
relate the class  $[\ib{\kappa}{\mu} ]$
with
$\ib{\kappa}{N}$.

 From 2.1 (a) and (b), we have that
 $\{M,N\backslash\overline{U}\}$ is a saturated open covering of $M$.
 It possesses  a
subordinated partition of the unity $\{f,g\}$ made up with basic
functions (cf. \cite{RSW}). Consider $\nu$ a D-metric on $N$, which always exists since
$N$ is compact.
So, the metric
$$
\lambda = f \mu + (1-f) \nu
$$
is a bundle-like metric on $N$ with $\lambda|_{U} = \mu|_{U}$, which
is a D-metric.
This gives
$$
\ib{\kappa}{\lambda,b}|_{U}
=
\ib{\kappa}{\mu|_{U}}
=
\ib{\kappa}{\mu} |_{U}
$$
(cf. (\ref{rest}) and  (\ref{bost})).
Denote by  $I \colon U \to M$ and $J \colon U \to N$ the natural inclusions.
We have

$$
I^{*}[\ib{\kappa}{\mu} ] = [\ib{\kappa}{\mu} |_{U}] = [\ib{\kappa}{\lambda,b}|_{U}]
=
J^{*}[\ib{\kappa}{\lambda,b}] = J^{*}\kappab_{N}.$$
Consider $\mu'$ another D-metric on $M$. The above equality gives
$
I^{*}[\ib{\kappa}{\mu} ]
=
I^{*}[\ib{\kappa}{\mu'}] .
$
From 2.1 (c) we get $
[\ib{\kappa}{\mu} ]
=
[\ib{\kappa}{\mu'}] .
$
\qed

The first characterization of the tautness is the following.
\bt
\label{T1}
Let $M$ be a manifold endowed with a CERF $\F$.
Then the following two statements are equivalent:

\Zati The foliation $\F$ is taut.

\zati The tautness class $\kappab \in \hiru{H}{1}{M/\F}$ vanishes.
\et
\pro
We prove the two implications.

\medskip

$(a) \Rightarrow (b)$. There exists a D-metric $\mu$ on $M$
with $\ib{\kappa}{\mu} = 0$. Then $\kappab = [\ib{\kappa}{\mu} ] = 0$.

\medskip

$(b) \Rightarrow (a)$. See \cite[Proposition
7.6]{To}\footnote{At the beginning of Chapter 7 of \cite{To} it is said
the the manifold $M$ must be compact. In fact, this condition is not
necessary on the proof of Proposition 7.6.}.

\qed

For the second characterization of the tautness we use the twisted
basic cohomology $\lau{H}{*}{\ib{\kappa}{\mu} }{\mf}$, where $\mu$ is a D-metric.
Notice that this cohomology does not depend on the choice of the
D-metric (cf.  Proposition  \ref{2D} and  (c) of page \pageref{zazpi}).

\bt
\label{T2}
Let $M$ be a manifold endowed with a CERF $\F$. Consider $\mu$ a
D-metric on $M$.
Then, the  following two  statements are equivalent:

\Zati The foliation $\F$ is taut.

\zati The cohomology group $\lau{H}{0}{\ib{\kappa}{\mu} }{M/\F}$ is $\R $.
\bigskip

\nt Otherwise, $\lau{H}{0}{\ib{\kappa}{\mu} }{\mf}= 0$.
\et
\pro
We proceed in two steps.

\medskip

$(a) \Rightarrow (b)$. If $\F$ is taut then $\kappab =
[\ib{\kappa}{\mu} ]=0$. So,
$\lau{H}{0}{\ib{\kappa}{\mu} }{\mf} \cong \hiru{H}{0}{\mf} = \R$.

\smallskip

$(b) \Rightarrow (a)$.  If $\lau{H}{0}{\ib{\kappa}{\mu} }{\mf}  \neq 0$
then there exists $0 \neq f \in \hiru{\Om}{0}{\mf}$ with $df = f
\ib{\kappa}{\mu} $. The set $Z(f) = f^{-1}(0)$  is clearly a closed subset of
$M$. Let us see that it is also an open subset. Take $x \in Z(f)$ and
consider a contractible open subset $V \subset M$ containing $x$. So,
there exists a smooth map $g \colon V \to \R$ with $\ib{\kappa}{\mu} =
dg $ on $U$. The calculation
$$
d(f e^{-g}) = e^{-g} df - fe^{-g} dg = e^{-g} f \ib{\kappa}{\mu} - e^{-g} f
\ib{\kappa}{\mu} =0
$$
shows that $f e^{-g}$ is constant on $V$. Since $x \in Z(f)$ then $f
\equiv 0$ on $V$ and therefore $x \in V \subset Z(f)$. We have proved  that
$Z(f)$ is an open subset.

As $M$ is connected, we have $Z(f) = \emptyset$. From
${\displaystyle d(\log |f|) = \frac{1}{f} df = \ib{\kappa}{\mu} }$
we conclude that $\kappab =0$. The foliation $\F$ is taut.

\medskip

Notice that we have also proved:  $\lau{H}{0}{\ib{\kappa}{\mu} }{\mf}  \neq 0
\Rightarrow \lau{H}{0}{\ib{\kappa}{\mu} }{\mf}  = \R$.
\qed

\prgg {\bf Remark.} To prove $(b) \Rightarrow (a)$ we do not need
$\F$ to be a CERF but just the existence of a $D$-metric $\mu$ (see
\cite[Proposition 7.6]{To}).

\medskip

For the third characterization we need to extend the basic
Poincaré duality to the non-compact case. We find in \cite{Se}
another version of this Poincaré duality using the cohomological
orientation sheaf instead of the twisted basic cohomology we use
here. Also compare with \cite[Proposition 7.54]{To}.

\bt
 \label{PK} Let $M$ be a manifold endowed with a transversally
oriented RF $\F$ possessing a zipper.
 Consider $\mu$ a $D$-metric on $M$.
If $n = \codim \F$, then
$$
\lau{H}{*}{c}{\mf} \cong\lau{H}{n - *}{\ib{\kappa}{\mu} }{\mf}.
$$
\et
\pro
See Appendix I.
   \qed

The third characterization of the tautness is the following. Compare with
\cite[Proposition 7.56]{To}.

\bt
\label{T3}
Let $M$ be a manifold endowed with a CERF $\F$.
Let us suppose that $\F$ is transversally oriented.
If $n = \codim \F$, then
the two following statements are equivalent:

\Zati The foliation $\F$ is taut.

\zati The cohomology group $\lau{H}{n}{c}{M/\F}$ is $\R $.

\bigskip

\nt Otherwise, $\lau{H}{n}{c}{\mf}= 0$.
\et
\pro It suffices to apply
Theorem \ref{PK} and Theorem
\ref{T2}.
\qed

As a direct application, we extend the scope of a well known
result for closed manifolds
(cf. \cite[Corollary 6.6]{Suso}) to arbitrary CERFs.

\bc Any codimension one CERF is taut.
\ec
 \pro Let $\F$ be a codimension one  CERF defined on a manifold $M$.
Without loss of generality we can suppose that $\F$ is transversally oriented (cf. \cite[Lemma 6.3]{Suso}).
By {\em reductio ad absurdum}, let's  suppose that
$\F$ is not taut. Then, by Theorem \ref{T2}, we get  $\lau{H}{0}{\ib{\kappa}{\mu} }{M/\F}= 0 $, and by
Theorem \ref{T1}, we get $\kappab \ne 0$, thus  $\lau{H}{1}{}{M/\F}\ne 0 $.
 Now, from   Remark 4.11 (b) we get $\lau{H}{0}{\ib{\kappa}{\mu}, c }{\mf} \ne 0$, and  then
$\lau{H}{0}{\ib{\kappa}{\mu} }{\mf} \ne 0$. Theorem \ref{PK} yields
$\lau{H}{1}{c}{\mf} \ne 0$, a contradiction.\qed

% Now, if
% $\F$ is not taut, then $\lau{H}{0}{\ib{\kappa}{\mu}}{M/\F}=0$.

%
\section{Appendix.}

This appendix is devoted to the  proof the Theorem \ref{PK}. We
distinguish two cases following the orientability of $M$.
Beforehand, we introduce two technical tools.

\prg {\bf Bredon's Trick.}
The Mayer-Vietoris sequence allows  us to make computations when the
manifold is covered by a finite suitable covering. The passage from
the finite case to the general case may be done using an adapted version of
the Bredon's trick of  \cite[page 289]{Br}  we present now.

Let $X$ be a paracompact topological space
 and  let $\{ U_\alpha\}$ be an open covering, closed for
 finite intersections. Suppose that $Q(U)$ is a statement about open subsets
 of $X$, satisfying the following three properties:

\begin{itemize}
    \item[]
    \begin{itemize}
    \item[(BT1)] $Q(U_\alpha)$ is true for each $\alpha$;

    \item[(BT2)] $Q(U)$, $Q(V)$ and $Q(U\cap V )$ $\Longrightarrow$ $Q(U\cup V
 )$, where $U$ and $V$ are open subsets of $X$;

 \item[(BT3)] $Q(U_i) \Longrightarrow Q\left({\displaystyle  \ \bigcup_i
 }U_i\right)$,
 where $\{ U_i\}$ is an arbitrary disjoint family  of open subsets of $X$.
\end{itemize}
\end{itemize}
 Then $Q(X)$ is true.

\medskip

\prg {\bf Mayer-Vietoris.}
Associated to a covering $\{ U , V \}$  of $M$ by saturated open
subsets we have the Mayer-Vietoris exact short sequence
$$
0 \to (\hiru{\Om}{*}{\mf} ,d)\to
(\hiru{\Om}{*}{U/\mathcal{F}},d) \oplus
(\hiru{\Om}{*}{V/\mathcal{F}},d)  \to
(\hiru{\Om}{*}{(U \cap V)/\mathcal{F}},d)  \to 0,
$$
where the maps are defined by restriction (see for example
\cite{RSW}). In the compact support context we have the
Mayer-Vietoris sequence
 $$
 0 \to  (\lau{\Om}{*}{c}{(U \cap V)/\mathcal{F}},d) \to
 (\lau{\Om}{*}{c}{U/\mathcal{F}},d) \oplus
 (\lau{\Om}{*}{c}{V/\mathcal{F}},d)  \to
 (\lau{\Om}{*}{c}{\mf} ,d)\to 0,
 $$
 where the maps are defined by extension
 (see for example \cite{RSW}).
 Finally, for the twisted basic cohomology, we have
  the Mayer-Vietoris sequence
 $$
 0 \to (\hiru{\Om}{*}{\mf} ,\ib{d}{\ib{\kappa}{\mu}})\to
 (\hiru{\Om}{*}{U/\mathcal{F}},\ib{d}{\ib{\kappa}{\mu}}) \oplus
 (\hiru{\Om}{*}{V/\mathcal{F}},\ib{d}{\ib{\kappa}{\mu}})  \to
 (\hiru{\Om}{*}{(U \cap V)/\mathcal{F}},\ib{d}{\ib{\kappa}{\mu}})  \to 0,
 $$
 where the maps are defined by restriction.

 \bigskip

\begin{center}
\bf
- - - Orientable case - - -
\end{center}

 \prg {\bf Integration.}
In order to define the duality operator, we fix

\Zati  an oriented manifold $M$,

\zati a transversally oriented RF $\F$ (TORF for short) on
$M$, and

\zati a $D$-metric $\mu$ on $(M,\F)$.

\smallskip

\nt We shall say that $(M,\F,\mu)$ is a {\em $D$-triple}.
 The  associated  tangent volume form  is $\ib{\chi}{\mu}$  (it exists since $\F$ is also oriented).  With all these
ingredients we  define the
morphism
$$
\Int{M}{}\colon \lau{H}{*}{c}{\mf}
\TO \Hom \left(\lau{H}{n-*}{\ib{\kappa}{\mu} }{\mf} ;\R\right)
=
\left( \lau{H}{n-*}{\ib{\kappa}{\mu} }{\mf}
\right)^{\star}
$$
by $\Int{M}{}([\alpha])([\beta]) = \Int{M}{} \alpha \wedge \beta
\wedge \ib{\chi}{\mu}$. Here $n = \codim_{M}\F$.
This operator
 is well defined since $M$ is oriented and we have the Rummler
formula
\begin{equation}
    \label{Rum}
    i_{Y_{1}} \cdots i_{Y_{r}}d \ib{\chi}{\mu} +\ib{\chi}{\mu}(Y_{1}, \ldots ,
    Y_{r})  \cdot \ib{\kappa}{\mu} = 0
    \end{equation}
     when  $\{Y_{1},
    \ldots , Y_{r}\} $ are vector fields tangent to  $T\F$ and $r =
    \dim \F$
    (see \cite{To}).
    We prove in this section that the operator $\Int{M}{}$ is an
    isomorphism.

\medskip

 Before proving the general case  first  we consider some particular
 cases.
     \bl
    \label{ER}
    Suppose that the $D$-triple $(M,\F,\mu) $ is $(E \times \R,\mathcal{E} \times
    \mathcal{I},\mu)$
    where  $E$ is a closed manifold
and the leaves of
 $\mathcal{E}$ are dense.
    Then the operator $\Int{M}{}$ is an isomorphism.
    \el
 \pro
For the proof of Lemma we proceed
 in several steps. We shall use the following notation.
 Given a differential form (or Riemannian metric) $\om$ on $E \times
 \R$ we shall denote by  $\om (t)$ the restriction $I_{t}^{*}\om$, where $I_{t}\colon E \to E
 \times \R$ is defined by $I_{t}(x) = (x,t)$ for each $t \in \R$.

 \medskip

 \underline{\em First Step}. {\em The cohomology
 $\lau{H}{*}{c}{(E \times \R)/\mathcal{E} \times \mathcal{I}} $}.
 Consider $f \colon \R \to [0,1]$ a smooth function with compact
 support such that $\Int{\R}{} f dt =1$. We know that the
correspondence $[\gamma] \mapsto [f \gamma  \wedge dt]$ establishes
an
 isomorphism between $\lau{H}{*}{}{E/\mathcal{E} } $
 and $\lau{H}{*+1}{c}{(E \times \R)/\mathcal{E} \times
 \mathcal{I}} $ (cf. \cite{RSW}).
 In fact, this isomorphism  does not depend on choice of $f$.

 \smallskip

 \underline{\em Second Step}. {\em The cohomology
 $ \lau{H}{*}{\ib{\kappa}{\mu} }{(E \times \R) / \mathcal{E} \times \mathcal{I}} $}.
 Notice that the metric $\mu(0)$ is a $D$-metric (cf. \refp{hamar}).
 So, we have two $D$-metrics on $E \times \R$: $\mu$ and $\mu (0) +
 dt^{2}$
 with $\ib{\kappa}{\mu(0) + dt^{2}} = \ib{\kappa}{\mu(0)}$.
 Since $\mathcal{E} \times \mathcal{I}$  is a CERF (it suffices
      to take $E \times ]-1,1[$ as a
    reppiz and $(E \times \S^{1},\mathcal{E} \times \I)$ as a zipper)
    then there
exists a function $g \in
 \hiru{\Om}{0}{(E \times \R) / \mathcal{E} \times \mathcal{I}}$ with
 $
 \ib{\kappa}{\mu} = \ib{\kappa}{\mu(0)} +  dg
 $
(cf.
  Proposition \ref{2D}).
Since the leaves of $\mathcal{E}$ are dense on $N$ then the (basic) function $g$  is a smooth function on $\R$.

 We know (see  (c) of page \pageref{zazpi}) that the assignment
 $[\om] \mapsto [e^g \om]$ establishes an isomorphism
 between
 $
 \lau{H}{*}{\ib{\kappa}{\mu(0)}}{(E \times \R) / \mathcal{E} \times
 \mathcal{I}}
$
and
$
 \lau{H}{*}{\ib{\kappa}{\mu} }{(E \times \R) / \mathcal{E} \times \mathcal{I}}
 $
 .
 The usual techniques give that the assignment $[\omega] \mapsto
 [\omega]$ establishes an isomorphism between
 $
 \lau{H}{*}{\ib{\kappa}{\mu(0)} }{E/ \mathcal{E} }
$
and
$
\lau{H}{*}{\ib{\kappa}{\mu(0)}  }{(E \times \R) / \mathcal{E} \times
\mathcal{I}}.
$

\smallskip

 \underline{\em Last Step}. Notice that
 $(E,\mathcal{E},\mu (0))$ is a $D$-triple.
Since $E$ is compact, the morphism
\begin{equation}
    \label{IntE}
  \Int{E}{} \colon \lau{H}{*}{}{E/\mathcal{E}} \TO
  \left( \lau{H}{n-1-*}{\ib{\kappa}{\mu(0)}}{E/\mathcal{E}}
  \right)^\star,
  \end{equation}
  defined by $\Int{E}{}([\gamma])([\zeta]) = \Int{E}{} \gamma
 \wedge \zeta
  \wedge \ib{\chi}{\mu}(0)$,
  is an isomorphism (see \cite{KT_AST}).
 Following the previous steps  it
 suffices to show  that the morphism
 $$
 \Int{E}{\#} \colon \lau{H}{*}{}{E/\mathcal{E}} \TO
 \left(\lau{H}{n-1-*}{\ib{\kappa}{\mu(0)}}{E/\mathcal{E}}
\right)^\star,
 $$
 defined by $\Int{E}{\#}([\gamma])([\zeta]) = \Int{E \times \R}{} f
 e^g\gamma
 \wedge dt \wedge \zeta
 \wedge \ib{\chi}{\mu}$,
 is an isomorphism. Let us see that.

 \begin{itemize}
     \item[-] $\Int{E}{\#}$ is a monomorphism.
     \end{itemize}
     Consider $[\gamma] \in \hiru{H}{*}{E/\mathcal{E}}$ with
     $\Int{E}{\#} ([\gamma])\equiv 0$. We have
     $\Int{\R}{}  f(t) e^{g(t)}\left( \Int{E}{}\gamma
 \wedge \zeta
 \wedge \ib{\chi}{\mu}(t)\right) dt  = 0
 $
 for each $[\zeta] \in
 \lau{H}{n-1-*}{\ib{\kappa}{\mu(0)}}{E/\mathcal{E}}$ and each smooth
 function $f \colon \R
 \to [0,1]$
 with compact support and $\Int{\R}{}f =1$. So, $\Int{E}{}\gamma\wedge \zeta
 \wedge \ib{\chi}{\mu(t)}dt  = 0$ for each $[\zeta] \in
 \lau{H}{n-1-*}{\ib{\kappa}{\mu(0)}}{E/\mathcal{E}}$ and each $t \in \R$.
 We get     $\Int{E}{}([\gamma])([\zeta ])
 = 0$ for each $[\zeta] \in
 \lau{H}{n-1-*}{\ib{\kappa}{\mu(0)}}{E/\mathcal{E}}$. Since
 $\Int{E}{}$ is an isomorphism then $[\gamma ] =0$.

 \begin{itemize}
 \item[-] $\Int{E}{\#}$ is an epimorphism.  From \refp{IntE} we know
 that
 $ \dim \lau{H}{*}{}{E/\mathcal{E}} =
 \dim \left(\lau{H}{n-1-*}{\ib{\kappa}{\mu(0)}}{E/\mathcal{E}}
\right)^\star,
$
which is finite since $E$ is compact (see \cite{EHS}). This gives that
the
monomorphism
$\Int{E}{\#}$ is also  an epimorphism.
\qed
 \end{itemize}

 \bl
 \label{TP}
 Let $(M,\mu,\F)$ be a $D$-triple.  Suppose that $(M,\F)$ possesses a
 transversally parallelizable zipper. Then $\Int{M}{}$ is an isomorphism.
 \el
 \pro
Denote by $(N,\HH)$ the transversally parallelizable zipper.
Since $\HH$ is transversally parallelizable, then there exists a fiber
 bundle $\pi \colon N \to B$ whose fibers are the closures of the
 leaves of $\HH$. Since $M$ is saturated for the leaves of $\HH_{M}$ then
 it is also saturated for the  closures of these leaves. We get an
 open subset $V_{M} \subset B$ with $M = \pi^{-1}(V_{M})$.

Fix $V$ an open subset of $V_{M}$. The foliation
 $\HH_{\pi^{-1}(V)}$
admits  the zipper $(N,\HH)$.
  Notice that $\pi^{-1}(V)$ is oriented and  $\HH_{\pi^{-1}(V)}$
  is a TORF.  The triple $(\pi^{-1}(V), \HH_{\pi^{-1}(V)},\mu_{\pi^{-1}(V)})$
is a $D$-triple. So, the operator $\Int{\pi^{-1}(V)}{}$ is
  well defined. We prove that this operator is non-degenerate. This
  will end the proof by taking $V = V_{M}$.

Let  $(E,\mathcal{E})$  be a generic fiber of $\pi$. The manifold
 $E$ is closed and the leaves of $\mathcal{E}$ are dense in $E$.
 We know that the fibration $\pi \colon \pi^{-1}(V) \to V$
 has a foliated atlas $\mathcal{A} =\{ \phii \colon
 (\pi^{-1}(U),\HH) \TO (U \times E,\mathcal{I} \times
 \mathcal{E}) \}$. We can suppose that the covering $\mathcal{U} = \{ U
 \ | \
 \exists (U,\phii) \in \mathcal{A}\}$:
 \begin{itemize}
 \item[-]  is a good covering of $V$: if $U_{1}, \ldots, U_{k} \in \mathcal{U}$ then the
 intersection $V = U_{1} \cap \cdots \cap U_{k} $ is diffeomorphic to $\R^{\dim B}$
 (cf.
  \cite{BT}), and
  \item[-] is closed for finite intersections.
  \end{itemize}
 We consider the statement $Q(U)$:
 \smallskip

 \begin{center}
     ``The integration operator
     $
    \Int{\pi^{-1}(U)}{} $
     is an isomorphism.''
     \end{center}

     \smallskip

 \nt where $U \subset V$ is an open subset. Following the Bredon's
 trick, it suffices to prove
 (BT1), (BT2) and (BT3) relatively to the covering $\mathcal{U}$.

 \medskip

 $+$ (BT1). It follows directly from  the Lemma \ref{ER}.

 \smallskip

 $+$ (BT2). The integration operator $\Int{}{}$ commutes with the
 restriction and inclusion operators. It suffices to apply the Five's Lemma to the
 Mayer-Vietoris sequences of 4.2.

 \smallskip

  $+$ (BT3). Straightforward.
 \qed

 \prg{\bf Frame bundle.}  Let $(M,\mu,\F)$ be a $D$-triple
 possessing a zipper $(N,\HH)$.
 Consider $p \colon \wt{N} \to N$ the bundle of transverse oriented
 orthonormal  frames  of $N$ (cf. \cite{Mol}). It is an
 $SO(n)$-principal bundle.
 The canonical lift $\wt{\HH}$ of $\HH$ is a transversally parallelizable
 foliation on the closed manifold $\wt{N}$ with $\codim_{\wt{N}}
 \HH =n + \dim SO(n)$.
 The  restriction bundle morphism $p_* \colon
 T \wt{\HH} \to T\HH$ is an isomorphism.
 We can lift $(M,\F,\mu)$ as follows.

 \medskip

- {\em Lifting $\F$.}
Since $M$ is a saturated open subset of $(N,\HH)$ then $\wt{M} = p^{-1}(M)$ is a saturated
open subset
 of $(\wt{N},\wt{\HH})$.  The foliation $\wt{\F} =\wt{\HH}_{\wt{M}}$
 is transversally parallelizable (and then a TORF)
 and the manifold $\wt{M}$ is oriented since $p \colon \wt{M} \to M$
 is a $SO(n)$-bundle and $M$ is oriented.
 The foliation $\wt\F$ possesses $(\wt{N},\wt\HH)$ as a zipper.

 \smallskip

-  {\em Lifting $\mu$}.
 Consider the decomposition $\mu = \mu_1 + \mu_2$ relatively to the
 orthogonal decomposition $TM = T\F \oplus \left( T\F
 \right)^{\ib{\bot}{\mu}}$.  Since the restriction  bundle morphism  $p_{*}\colon T \wt\F
 \to T\F$ is an isomorphism  then we have the decomposition
 $T\wt{M} = T\wt\F \oplus p_{*}^{-1}
 \left( T\F
  \right)^{\ib{\bot}{\mu}}$.
  Moreover, since $(\wt{M}, \wt\F)$ is a Riemannian foliated
  manifold (TP in fact) then there exists a Riemannian metric
  $\nu_2$ on $p_{*}^{-1}
 \left( T\F
  \right)^{\ib{\bot}{\mu}}$ such that the Riemannian metric $\nu =
  p^*\mu_{1} + \nu_{2}$ is a bundle-like metric on $(\wt{M},
  \wt\F)$. Then, the associated volume forms verify:
  %\begin{equation}
   %   \label{16}
    $$
      \ib{\chi}{\nu} = p^*\ib{\chi}{\mu}.
    $$
    %  \end{equation}
 Rummler's formula \refp{Rum} gives
 \begin{equation}
     \label{17}
 \ib{\kappa}{\nu} = p^{*}\ib{\kappa}{\mu} .
\end{equation}

\smallskip

We conclude that $(\wt{M},\wt\F,\nu)$ is a $D$-triple possessing a
transversally parallelizable zipper. From the Lemma \ref{TP} we know
that the integration operator $\Int{\wt{M}}{} \colon
\lau{H}{*}{c}{\wt{M}/\wt\F} \TO \left( \lau{H}{n+
\ell-*}{\ib{\kappa}{\nu} }{\wt{M}/\wt\F} \right)^{\star} $ is an
isomorphism. Here $\ell = \dim SO(n)$. We shall prove Theorem
\ref{PK} relating $(M,\F,\mu)$ with $(\wt{M},\wt\F,\nu)$  using
two spectral sequences.

\prg{\bf A spectral sequence.}\footnote{This is the spectral
sequence of \cite[9.1,Ch.IX,vol.III]{GHV}.} Consider the usual
filtration
$$
 \bi{F}{p}
\lau{\Om}{p+q}{c}{\wt{M}/\wt{\F}}
 = \left\{ \om \in  \lau{\Om}{p+q}{c}{\wt{M}/\wt{\F}}
  \ \big/ \
\ib{i}{X_{u_0}} \cdots \ib{i}{X_{u_q}}\om
 = 0 \hbox{
 for each  } \{ u_{0}, \ldots , u_{q}\} \subset \son
 \right\},
$$
where   $X_{u}\in \mathfrak{X}(\wt{M})$ is  determined  by the
element $u\in \son$, the Lie algebra of $SO(n)$. It induces a
filtration in the differential complex
 $
 \KK^{^{*}}  =\left(   \lau{\Om}{*}{c}{\wt{M}/\wt{\F}}
 \right)^{SO(n)}
 $
by
$$
\bi{F}{p} \KK^{^{p+q}}= \KK^{^{p+q}} \cap \bi{F}{p}
\lau{\Om}{*}{c}{\wt{M}/\wt{\F}} ,
 $$
leading us to  a first quadrant spectral sequence $
\tres{\EE}{p,q}{r}$ which verifies
\begin{itemize}
    \item[(a)]  $
\tres{\EE}{p,q}{r} \Rightarrow \lau{H}{p+q}{c}{\wt{M}/\wt{\F}}$,
and

\item[(b)] $ \tres{\EE}{p,q}{2} \cong   \lau{H}{p}{c}{\mf}
\otimes
  \hiru{H}{q}{SO(n)}$.
\end{itemize}

Let us see that. The inclusion
 $
 \KK^{^{*}}   \hookrightarrow \lau{\Om}{*}{c}{\wt{M}/\wt{\F}}
 $
 induces an isomorphism in cohomology. This is a standard argument
 based on the fact that $SO(n)$ is a connected compact Lie group
 (cf. \cite[Theorem I,Ch.IV,vol.II]{GHV}). This gives (a).

 We denote by $\ib{\gamma}{u} = \ib{i}{X_{u}}\nu$
the associated fundamental differential form. Notice that the
 assignment $\alpha \otimes u \mapsto \alpha
 \wedge \gamma_u$
 induces the identification
 \begin{equation}
     \label{p1}
\bigoplus_{p+q=*} \left(
\bi{F}{p} \lau{\Om}{p}{c}{\wt{M}/\wt{\F}}
\otimes
{\bigwedge}^{q}\son
  \right)
= \lau{\Om}{*}{c}{\wt{M}/\wt{\F}}.
\end{equation}
  So, we get
$ \tres{\EE}{p,q}{0} \cong \left(\bi{F}{p}
\lau{\Om}{p}{c}{\wt{M}/\wt{\F}} \otimes
  {\bigwedge}^{q}\son\right)^{SO(n)}
 .
$ A straightforward calculation gives (see also
\cite[(9.2),vol.III]{GHV}) that $d_{0} = -  \Ide \otimes \
\delta$, where $\ib{d}{0}$ is the $0$-differential of the spectral
sequence and $\delta$ is the differential of
${\bigwedge}^{*}\son$. This gives $ \tres{\EE}{p,q}{1} \cong
\left( \bi{F}{p} \lau{\Om}{p}{c}{\wt{M}/\wt{\F}}
\right)^{SO(n)}\otimes
  \lau{H}{q}{}{SO(n)}
$
 (cf. \cite[5.28 and 5.12, vol.III]{GHV}). On the other hand, we have
 $$
 \left( \bi{F}{p} \lau{\Om}{p}{c}{\wt{M}/\wt{\F}} \right)^{SO(n)}
 =
 \left\{ \om\in   \lau{\Om}{p}{c}{\wt{M}/\wt{\F}}
 \ \big/ \ \ib{i}{X_{u}} \om = \ib{L}{X_{u}} \om =0 \hbox{
 for each } u \in  \son \right\}
 =
 p^{*}\lau{\Om}{p}{c}{{M}/{\F}}
$$
and then $ \tres{\EE}{p,q}{1} \cong\lau{\Om}{p}{c}{{M}/{\F}}
\otimes
  \lau{H}{q}{}{SO(n)}
$. Since $\ib{d}{1}$, the $1$-differential of the spectral
sequence, becomes $d \otimes \Ide$ then we conclude
$\tres{\EE}{p,q}{2} \cong   \lau{H}{p}{c}{\mf}  \otimes
  \hiru{H}{q}{SO(n)}$. This gives (b).

  \prg{\bf Another spectral sequence.}\footnote{This is the spectral sequence
  of \cite[Ch.,IX,9.1,vol.III]{GHV} associated to
  $\left( \lau{\Om}{*}{\ib{\kappa}{\nu}}{\wt{M}/\wt{\F}}\right)^\star $
  which is \underline{not} a differential graded algebra.}
  Consider the usual
filtration
$$
 \bi{F}{p}
\lau{\Om}{p+q}{}{\wt{M}/\wt{\F}}
 = \left\{ \om \in  \lau{\Om}{p+q}{}{\wt{M}/\wt{\F}}
  \ \big/ \
\ib{i}{X_{u_0}} \cdots \ib{i}{X_{u_q}}\om
 = 0 \hbox{
 for each  } \{ u_{0}, \ldots , u_{q}\} \subset \son
 \right\}.
$$
 The complex $\KKK^{^{*}} =
\left( \left( \lau{\Om}{n+\ell -*}{}{\wt{M}/\wt{\F}}
   \right)^{SO(n)} \right)^\star
   $
   admits
 $\ib{\nabla}{\ib{\kappa}{\nu}}$, the dual of
$\ib{d}{\ib{\kappa}{\nu}}$, as a differential since \refp{17}.
Consider the filtration
$$
\bi{F}{p} \KKK^{^{p+q}}= \left\{ L \in \KKK^{^{p+q}} \ / \ L
\equiv 0 \hbox{ on } \left( \bi{F}{n-p+1}\lau{\Om}{n+\ell
-(p+q)}{}{\wt{M}/\wt{\F}}\right)^{SO(n)} \right\}.
$$
 A straightforward calculation gives that $\bi{F}{p+1}
\KKK^{^{*}}\subset \bi{F}{p} \KKK^{^{*}}$ and
$\ib{\nabla}{\ib{\kappa}{\nu}}(\bi{F}{p} \KKK^{^{p+q}})  \subset
\bi{F}{p} \KKK^{^{p+q+1}} $.
 Thus, it induces a first quadrant spectral sequence
  $
  \tres{\EEE}{p,q}{r}$
  verifying
  \begin{itemize}
      \item[(a)]  $
  \tres{\EEE}{p,q}{r} \Rightarrow
 \left( \lau{H}{n+ \ell -(p+q)}{\ib{\kappa}{\nu}}{\wt{M}/\wt{\F}}\right)^\star.$

  \item[(b)] $
  \tres{\EEE}{p,q}{2} \cong   \left(\lau{H}{n-p}{\ib{\kappa}{\mu}}{\mf} \right)^\star \otimes
    \left(\hiru{H}{\ell -q}{SO(n)}\right)^\star$.
  \end{itemize}
  Let us see that. As in 4.7 (a),
  the inclusion
 $
\left( \lau{\Om}{*}{}{\wt{M}/\wt{\F}}
   \right)^{SO(n)}   \hookrightarrow  \lau{\Om}{*}{}{\wt{M}/\wt{\F}}
 $
 induces an isomorphism in the corresponding twisted basic cohomology.  This yields (a).

 The analogous identification to \refp{p1} is here
  %\begin{equation*}
   %  \label{p2}
$$
\bigoplus_{p+q=*} \left( \bi{F}{p} \lau{\Om}{p}{}{\wt{M}/\wt{\F}}
\otimes
  {\bigwedge}^{q}\son
  \right)
= \lau{\Om}{*}{}{\wt{M}/\wt{\F}}.
$$
%\end{equation*}
 So, we get
  $
  \tres{\EEE}{p,q}{0}
 \cong
 \left( \left(\bi{F}{n-p} \lau{\Om}{n-p}{}{\wt{M}/\wt{\F}} \otimes
    {\bigwedge}^{\ell -q}\son\right)^{SO(n)}\right)^\star
   .
  $
  A straightforward calculation shows
 % (see also \cite[(9.2),vol.III]{GHV})
  that
  the $0$-differential of the spectral
sequence is the dual of
  $-\Ide \otimes \ \delta$. This gives
  $
  \tres{\EE}{p,q}{1}
  \cong
  \left(\left( \bi{F}{n-p} \lau{\Om}{n-p}{}{\wt{M}/\wt{\F}} \right)^{SO(n)}\right)^\star\otimes
    \left(\lau{H}{\ell -q}{}{SO(n)}\right)^\star
  $
   %(cf. \cite[5.28 and 5.12,vol.III]{GHV})
   . On the other hand, we have
   $$
 \left( \bi{F}{n-p} \lau{\Om}{n-p}{}{\wt{M}/\wt{\F}} \right)^{SO(n)}
   =
   \left\{ \om\in   \lau{\Om}{n-p}{}{\wt{M}/\wt{\F}}
   \ \big/ \ \ib{i}{X_{u}} \om = \ib{L}{X_{u}} \om =0 \hbox{
   for each } u \in  \son \right\}
   =
   p^{*}\lau{\Om}{n-p}{}{{M}/{\F}}
  $$
  and then
  $
  \tres{\EE}{p,q}{1}
 \cong \left(\lau{\Om}{n-p}{}{{M}/{\F}} \right)^\star \otimes
    \left( \lau{H}{\ell - q}{}{SO(n)}\right)^\star
  $.
  Since  the
  $1$-differential of the spectral sequence becomes the dual of $\ib{d}{\ib{\kappa}{\mu}}\otimes
  \Ide$ (cf. \refp{17})
  then we conclude
  $\tres{\EE}{p,q}{2} \cong   \left(\lau{H}{n-p}{\ib{\kappa}{\mu}}{\mf}\right)^\star  \otimes
  \left(  \hiru{H}{\ell - q}{SO(n)}\right)^\star$.

 \bp
 \label{orient}
 Let $(M,\mu,\F)$ be a $D$-triple.
 Suppose that $\F$ possesses a zipper, then
 the operator $\Int{M}{}$ is an isomorphism.
 \ep
 \pro
Consider the differential operator
  $$\Delta \colon
  \left(\left( \lau{\Om}{*}{c}{\wt{M}/\wt{\F} }\right)^{SO(n)};d\right)\TO
   \left(\left(  \left(
   \lau{\Om}{n + \ell - *}{}{\wt{M}/\wt{\F}}\right)^{SO(n)}
   \right)^\star;\ib{\nabla}{\ib{\kappa}{\nu}}\right)
   $$
 defined by
  $\Delta (\om)(\eta) = \Int{\wt{M}}{} \om \wedge \eta
 \wedge \ib{\chi}{\ib{\kappa}{\nu}}$ (cf. 4.3).  By degree reasons, it preserves the involved
 filtrations, that is, we have
 $\Delta \left(F^{p}\KK^{p+q}\right) \subset F^{p}\KKK^{p+q}$.
It induces the following morphisms:
\begin{itemize}
 \item[+] [At $\infty$-level]  $\Int{\wt{M}}{} \colon
  \lau{H}{p+q}{c}{\wt{M}/\wt{\F} }\TO
  \left(
   \lau{H}{n + \ell - (p+q)}{\ib{\kappa}{\nu}}{\wt{M}/\wt{\F}}
   \right)^\star$.
 \item[+] [At $2$-level] $\Int{M}{} \otimes \Int{SO(n)} {} \colon
 \lau{H}{p}{c}{\mf}  \otimes
  \hiru{H}{q}{SO(n)}\TO \left( \lau{H}{n
  -p}{\ib{\kappa}{\mu} }{\mf}\right)^\star
     \otimes \left(\hiru{H}{\ell - q}{SO(n)}\right)^\star$.
 \end{itemize}
 As the operators $\Int{\wt{M}}{}$  and $\Int{SO(n)}{}$
    are isomorphisms (cf. Lemma \ref{TP}),
    Zeeman's
comparison theorem yields that $\Int{M}{} $ is an isomorphism (see
for example \cite{McC}). \qed

\bigskip
\begin{center}
\bf
- - - Non-orientable case - - -
\end{center}

For the non-orientable case  it suffices to consider the orientation
covering in order to apply the previous results.

\bp
   Let $\F$ be a TORF defined on
   a  manifold $M$.  Suppose that $\F$ possesses a zipper.
 Consider $\mu$ a $D$-metric on $M$.
 If $n = \codim \F$, then
   $$
   \lau{H}{*}{c}{\mf} \cong\lau{H}{n - *}{\ib{\kappa}{\mu} }{\mf}.
   $$
   \ep
 \pro
Let us suppose that $M$ is not orientable. We fix  a  a zipper
$(N,\HH)$ of $\F$.
 We consider $\mu$ a $D$-metric on $(M,\F)$.

 Consider $\oslash \colon \check{N}
    \to N$ the two-fold orientation
      covering of $N$. It is an oriented closed manifold. Denote by  $\check{\HH}$ the lifted foliation,
      which is a Riemannian one. In fact, there exists a smooth foliated
      action $\Phi \colon \Z_{2} \times (\check{N},\check{\HH})\to (\check{N},\check{\HH})$ such
      that $\oslash$ is $\Z_{2}$-invariant and
      $\check{N}/\Z_{2} = N$. Put $\flat \colon (\check{N},\check{\HH})\to
      (\check{N},\check{\HH})$ the foliated diffeomorphism generating
      this action.

      The restriction $\oslash \colon \oslash^{-1} (M) \to M$ is the the two-fold orientation
      covering of $M$. The manifold $\check{M} = \oslash^{-1} (M)$ is oriented and
      $\check{\HH}$-saturated. The diffeomorphism $\flat \colon \check{M}
      \to \check{M}$ preserves the foliation $\check{\F}=
      \check{\HH}|_{\oslash^{-1} (M)}$ and the
      $D$-metric $\check{\mu} = \oslash^{*}\mu$.

      Since the foliation $\F$ is
      transversally oriented then the foliation $\check{\F}$ is
        also transversally oriented. Moreover, the diffeomorphism  $\flat
     \colon \check{M} \to \check{M}$ preserves the
      transversally orientation of $\check{\F}$. It does not preserve the
      orientation of $\check{M}$
      since $M$ is a not orientable manifold. We get that $\flat$ does
      not preserve the tangential orientation of $\check{\F}$, this
      gives:
      \begin{equation}
      \label{ji}
      \flat^{*}\ib{\chi}{\check{\mu}} = - \ib{\chi}{\check{\mu}}.
      \end{equation}

      The foliated manifold
      $\check{\F}$ is a TORF on an oriented manifold $\check{M}$ with
      $(\check{N},\check{\HH})$ as a zipper. So,
      $(\check{M},\check{\F},\check{\mu})$ is a $D$-triple.
      From Theorem
      \ref{PK} we get that $\Int{\check{M}}{}$ induces an
      isomorphism
      \begin{equation}
      \label{iso2}
      \lau{H}{*}{c}{\check{M}/\check{\F}} \cong
      \lau{H}{n - *}{\ib{\kappa}{\check{\mu}}}{\check{M}/\check{\F}} ,
      \end{equation}
      where $n = \codim_{\check{M}}\check{\F} = \codim_{M}\F$.

 On the other hand, the map $\oslash$ induces the isomorphisms
 $
 \lau{H}{*}{c}{\mf} \cong \left(
 \lau{H}{*}{c}{\check{M}/\check{\F}}\right)^{\Z_{2}}
 $
 and
 $
    \lau{H}{*}{\ib{\kappa}{\mu} }{\mf} \cong \left(
    \lau{H}{*}{\ib{\kappa}{\check{\mu}}}{\check{M}/\check{\F}}
    \right)^{\Z_{2}}
    $
    since $\flat^{*}\ib{\kappa}{\check{\mu}} =
    \ib{\kappa}{\check{\mu}}$.
    From \refp{iso2} it suffices to prove that $\Int{\check{M}}{}$
    is $\Z_{2}$-invariant. This comes from the equality
    $$
    \Int{\check{M}}{} \flat^{*} \alpha \wedge \beta \wedge \ib{\chi}{\check{\mu}}
    \stackrel{\refp{ji}}{=\! = }
    - \Int{\check{M}}{} \flat^{*} \alpha \wedge \beta \wedge \flat^{*} \ib{\chi}{\check{\mu}}
    \stackrel{\flat^{-1} = \flat}{=\! = \! = \! =}
    -\Int{\check{M}}{} \flat^{*} \alpha  \wedge \flat^{*}\flat^{*}\beta \wedge \flat^{*} \ib{\chi}{\check{\mu}}
    \stackrel{\flat \footnotesize{\ not \ orient.}}{=\! = \! = \!
    =\!=\! =  \! =}
    \Int{\check{M}}{} \alpha \wedge \flat^{*}\beta \wedge
    \ib{\chi}{\check{\mu}} ,
    $$
    where $\alpha \in
    \lau{\Om}{*}{c}{\check{M}/\check{\F}}$ and
    $\beta \in \lau{\Om}{n -
    *}{}{\check{M}/\check{\F}}.$
 \qed

\prg {\bf Remarks}.

\Zati
The above proof gives also that the pairing
$
I \colon  \lau{H}{*}{c}{\mf} \oplus \lau{H}{n - *}{\ib{\kappa}{\mu} }{\mf}
\TO \R
$
defined by
$
I( [\alpha], [\beta]) = \Int{\check{M}}{} \oslash^{*}\alpha \wedge
\oslash^{*}\beta \wedge \ib{\chi}{\check{\mu}},
$
is non-degenerate.

\zati Under the assumptions of Theorem  \ref{PK} we also have $$
\lau{H}{*}{}{\mf} \cong\lau{H}{n - *}{\ib{\kappa}{\mu},c }{\mf},
$$
where the twisted cohomology is with compact supports.

\end{document}